\documentclass{amsart}
\usepackage{amssymb,amscd}

\usepackage{amsmath}
\usepackage{amsthm}

\newtheorem{theorem}{Theorem}
\newtheorem{lemma}[theorem]{Lemma}

\newtheorem{proposition}[theorem]{Proposition}
\numberwithin{theorem}{section}

\theoremstyle{definition}

\theoremstyle{remark}
\newtheorem{remark}[theorem]{Remark}

\numberwithin{equation}{section}

\newfont{\germ}{eufm10}

\newcommand\B{{\mathcal B}}

\def\ol{\overline}
\def\Z{\mathbb Z}

\def\la{\lambda}
\def\ka{\kappa}
\def\vphi{\varphi}

\def\veps{\varepsilon}

\def\kab{\overline{\ka}}
\def\lab{\overline{\la}}
\def\xb{\overline{x}}
\def\yb{\overline{y}}

\begin{document}

\title{Tropical $R$ and tau functions}

\author{A. Kuniba}
\address{Institute of Physics, University of Tokyo, Tokyo 153-8902, Japan}
\email{atsuo@gokutan.c.u-tokyo.ac.jp}

\author{M. Okado}
\address{Department of Informatics and Mathematical Science,
    Graduate School of Engineering Science,
Osaka University, Osaka 560-8531, Japan}
\email{okado@sigmath.es.osaka-u.ac.jp}

\author{T. Takagi}
\address{Department of Applied Physics, National Defense Academy,
Kanagawa 239-8686, Japan}
\email{takagi@nda.ac.jp}

\author{Y. Yamada}
\address{Department of Mathematics, Faculty of Science,
Kobe University, Hyogo 657-8501, Japan}
\email{yamaday@math.kobe-u.ac.jp}

\begin{abstract}
Tropical $R$ is the birational map that 
intertwines products of geometric crystals and satisfies the 
Yang-Baxter equation.
We show that the $D^{(1)}_n$ tropical $R$ 
introduced by the authors and its reduction to 
$A^{(2)}_{2n-1}$ and $C^{(1)}_n$ 
are equivalent to a system of bilinear 
difference equations of Hirota type.
Associated tropical vertex models admit solutions 
in terms of tau functions of the BKP and DKP hierarchies.
\end{abstract}


\maketitle

\section{Introduction}\label{sec:intro}

Let $\B = \{x = (x_1,\ldots, x_n) \}$ be
the set of variables.
The tropical $R$ for $A^{(1)}_{n-1}$ \cite{HHIKTT,Ki,Y} is the birational map 
$R: \B \times \B \rightarrow \B \times \B$ specified by
\begin{equation}\label{eq:RPP}
\begin{split}
&R(x,y) = (x',y'),\qquad x'_i=y_i \frac{P_i(x,y)}{P_{i-1}(x,y)}, \quad
y'_i=x_i \frac{P_{i-1}(x,y)}{P_i(x,y)},\\
&P_i(x,y)=\sum_{k=1}^{n} \prod_{j=k}^n x_{i+j} \prod_{j=1}^k y_{i+j},
\end{split}
\end{equation}
where all the indices are considered to be in $\Z/n\Z$.
It satisfies the inversion relation $R^2=id$ on $\B \times \B$ and the 
Yang-Baxter equation $R_1R_2R_1 = R_2R_1R_2$ on $\B \times \B \times \B$,
where $R_1(x,y,z) = (R(x,y),z)$ and 
$R_2(x,y,z) = (x,R(y,z))$. 
Given $(x,y)$, it is characterized 
as the unique solution to the equations on $(x',y')$:
\begin{equation}\label{eq:toda}
x_iy_i = x'_iy'_i,\qquad 
\frac{1}{x_i}+\frac{1}{y_{i+1}} = 
\frac{1}{x'_i}+\frac{1}{y'_{i+1}}
\end{equation}
with an extra constraint
$\prod_{i=1}^n(x_i/y'_i) = 
\prod_{i=1}^n(y_i/x'_i)= 1$.

A representation theoretical background for the tropical $R$ is 
provided by the geometric crystals \cite{BK} and their 
natural extrapolation into the affine setting 
as demonstrated in section 1 in \cite{KOTY}.
To recall it consider the matrix realization, where 
each element $x \in \B$ is associated with the matrix 
(see also \cite{KNY})
\begin{equation*}
M(x,\zeta) = \begin{pmatrix}
x^{-1}_1 &  & & -\zeta \\
-1 & x^{-1}_2 & & \\
& -1 & \ddots & \\
& & -1 & x^{-1}_n
\end{pmatrix}^{-1}
\end{equation*}
involving the spectral parameter $\zeta$.
Then the structure of $A^{(1)}_{n-1}$ geometric crystal 
on $\B$ is realized as simple matrix operations.
For instance the geometric Kashiwara operator
$e^c_i:(x_1,\ldots, x_n) 
\mapsto (\ldots, x_{i-1},cx_i,c^{-1}x_{i+1}, x_{i+2},\ldots)$ 
is induced by a multiplication of $M$ with certain unipotent matrices.
The product of matrices 
$M(x^{(1)},\zeta) \cdots M(x^{(j)},\zeta)$ 
corresponds to the product of geometric crystals 
$(x^{(1)},\ldots, x^{(j)}) \in 
\B\times \cdots \times \B$ 
in the sense of \cite{BK}.
Eq.(\ref{eq:toda}) is equivalent to the matrix equation
\begin{equation}\label{eq:MM}
M(x,\zeta)M(y,\zeta) = M(x',\zeta)M(y',\zeta).
\end{equation}
Due to the presence of the spectral parameter $\zeta$, 
its non-trivial solution is unique, 
which characterizes the tropical $R$ as the intertwiner 
of the products of geometric crystals.

Geometric crystals are so designed that 
the relevant rational functions become totally positive 
\cite{BK,BFZ,L}.
There is no minus sign in 
$P_i(x,y)$ in (\ref{eq:RPP}) indeed.
Such functions can consistently be transformed into 
piecewise linear ones by replacing $+,\times$ and $/$ 
by $\max(\min), +$ and $-$, respectively.
Consequently, the structure of geometric crystals 
reduces to the one for crystal bases \cite{K}.
In the present case, the expression (\ref{eq:RPP})
leads to the piecewise linear formula 
\cite{HHIKTT} for the combinatorial $R$ \cite{NY}.
An analogous result is available in \cite{KOTY} 
for the $D^{(1)}_n$ combinatorial $R$ \cite{HKOT}.
The above transformation can be realized as a certain 
limiting procedure \cite{TTMS} and is often called 
the ultradiscretization (UD).
As related topics, we remark the 
tropical combinatorics \cite{Ki,NoY} and the soliton cellular automata 
associated with crystal bases \cite{TS,TNS,HKT1,FOY,HHIKTT,HKOTY,HKT2}.

Now let us turn to the aspect of the tropical $R$
as a classical integrable system, 
which is the main subject of this paper.
Recall the discrete time Toda equation \cite{HT,HTI} on the 
electric ``current" $I^t_i$ and ``voltage" $V^t_i$:
\begin{equation*}
I^{t+1}_iV^{t+1}_i = I^t_{i+1}V^t_i,\qquad 
I^{t+1}_i + V^{t+1}_{i-1} = 
I^t_i + V^t_i
\end{equation*}
with periodicity $I^t_i = I^t_{i+n}$, 
$V^t_i = V^t_{i+n}$.
This difference equation is known to be integrable.
Eq.(\ref{eq:toda}) is identified with it 
via $x^{-1}_i = I^t_{i+1}$,
$y^{-1}_i = V^t_i$, ${x'}^{-1}_i = V^{t+1}_i$ and 
${y'}^{-1}_i = I^{t+1}_i$ \cite{Y}.
The matrix equation (\ref{eq:MM}) is a Lax representation in this context.
In fact the tropical $R$ is equivalent to a system of 
bilinear difference equations of Hirota type, which was effectively 
the base of the analyses in \cite{HHIKTT}.
To see this, introduce the functions 
$\tau^J_i \, (1 \le J \le 4, i \in \Z/n\Z)$ and 
the parameters $\la_i, \ka_i$, and make the change of variables
\begin{align*}
&x^{-1}_i = \la_i\delta\tau^3_i/\delta\tau^2_i,\quad
y^{-1}_i = \ka_i\delta\tau^2_i/\delta\tau^1_i,\\
&{x'}^{-1}_i = \ka_i\delta\tau^3_i/\delta\tau^4_i,\quad
{y'}^{-1}_i = \la_i\delta\tau^4_i/\delta\tau^1_i
\end{align*}
with $\delta\tau^J_i = \tau^J_i/\tau^J_{i-1}$.
Then the former relation in (\ref{eq:toda}) is automatically satisfied
and the latter is translated into 
\begin{equation}\label{eq:abl}
\la_i\tau^2_{i-1}\tau^4_i - \ka_i\tau^2_i\tau^4_{i-1} 
= \alpha\tau^1_i\tau^3_{i-1}
\end{equation}
for any nonzero parameter $\alpha$ independent of $i$.
The birational map $R:(x,y) \mapsto (x',y')$ is induced by 
an automorphism  $\tau^2_i \leftrightarrow \tau^4_i$, 
$\la_i \leftrightarrow \ka_i$, $\alpha \rightarrow -\alpha$
of (\ref{eq:abl}).
Eq. (\ref{eq:abl}) is a version of so-called 
Hirota-Miwa equation \cite{H,M}, 
a prototype discrete soliton equation to which the
well-developed machinery of 
free fermions and infinite dimensional Lie algebras \cite{DJM,JM}
can be applied.
The solutions are provided by tau functions of the KP hierarchy 
with a certain reduction and discretization of time variables.

In this paper we elucidate the 
classical integrability of such sort for  
the $D^{(1)}_n$ tropical $R$ introduced in \cite{KOTY} 
together with its reductions to 
$A^{(2)}_{2n-1}$ and $C^{(1)}_n$.
We show that they are equivalent 
to a quartet of $A$ type bilinear equations like (\ref{eq:abl}), 
and construct solutions in terms of tau functions in
soliton theory \cite{JM}.
The cases $D^{(1)}_n$ and $A^{(2)}_{2n-1}$ are related to 
reductions of the DKP hierarchy associated to the algebras 
$D_\infty$ and $D'_\infty$. 
A key is the relation between 
two kinds of tau functions 
that originate in the isomorphism 
$D_\infty \simeq D'_\infty$.
See Lemmas \ref{lem:jm7.6} and \ref{lem:sapix}.
The $C^{(1)}_n$ case corresponds to a further reduction 
to the BKP hierarchy associated with the algebras 
$B_\infty \simeq B'_\infty$.

We shall formulate tropical vertex models, where 
the tropical $R$ plays a role of local time evolution.
Namely, it is the two dimensional system on a square lattice 
where each edge is assigned an element of a geometric crystal
in such a way that those four surrounding a vertex are related by the
tropical $R$.
Our solutions to the bilinear equations 
are naturally extended to the tropical vertex models 
by duplicating the original lattice and 
attaching a Clifford group element to each face of it 
under a certain rule.
The Yang-Baxter equation for the tropical $R$ is 
naturally understood {}from such a 
point of view. 
See Remark \ref{rm:ybe}.

We expect that the results in this paper are 
glimpses of deeper relations yet to be explored 
between geometric crystals, crystal bases, 
solvable lattice models, 
discrete and ultradiscrete solitons and so forth.
A rough picture at present looks as follows.
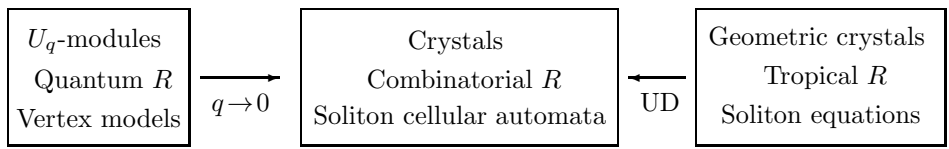
\begin{figure}[h]
\setlength{\unitlength}{1.5mm}
\begin{picture}(150,12)(-1,0)
\put(0,0){\framebox(16,12){%
\shortstack{$U_q$-modules \vspace{0.1cm} \\ \vspace{0.1cm}
Quantum $R$ \\ 
Vertex models}}}
\put(17,6){\vector(1,0){7}}
\put(18,3){$q\!\rightarrow\!0$}
\put(26,0){\framebox(28,12){%
\shortstack{Crystals \vspace{0.1cm}\\ \vspace{0.13cm}
Combinatorial $R$ \\ 
Soliton cellular automata}}}
\put(60,6){\vector(-1,0){5}}
\put(56,3){UD}
\put(61,0){\framebox(22,12){%
\shortstack{Geometric crystals \vspace{0.1cm} \\ \vspace{0.1cm}
Tropical $R$ \\ 
Soliton equations}}}
\end{picture}
\end{figure}

In \cite{D}, Drinfeld proposed 
set-theoretical solutions of the Yang-Baxter 
equation as one of 
unsolved problems in quantum group theory.
Since then a number of approaches have been made.
See for example \cite{ABS,ESS,JMY,MV,O,V1,V2,WX} and 
references therein.
Our tropical $R$ is a distinguishable example 
which is placed in the unique spot in the above picture and 
enjoys the total positivity.
It should be noted that the combinatorial $R$ is 
also a remarkable example of the set-theoretical solution 
which intertwines products of finite sets.

The paper is organized as follows.
In Section \ref{sec:TropicalR} we recall the 
tropical $R$ for $D^{(1)}_n$ following \cite{KOTY}.
In Section \ref{sec:BLE} the bilinear equations 
are introduced, which are divided into four groups 
corresponding to each 
vertex in Figure \ref{fig:fig1}.
Uniqueness and existence of the solution 
are shown in Proposition \ref{pr:unique-existence}. 
In Section \ref{sec:B-tropR} we establish the bilinearization of the 
tropical $R$ in Theorem \ref{th:bilinear}.
In Section \ref{sec:sol} we construct solutions of the 
bilinear equations in terms of tau functions 
of the DKP hierarchy formulated with 
two component fermions.
The final result is given in Theorem \ref{th:n-finite}.
In Section \ref{sec:red} we introduce the 
tropical $R$ for $A^{(2)}_{2n-1}$ and $C^{(1)}_n$ as 
reductions of the $D^{(1)}_n$ case. 
Parallel results on the bilinearization and 
solutions are obtained.
Appendix \ref{sec:appA} provides a proof of Lemma \ref{lem:UV}.
Appendix \ref{sec:appB} contains 
elements of the free fermion approach 
\cite{JM} as well as the key lemmas \ref{lem:jm7.6} and \ref{lem:sapix}.
We leave the calculation of the ultradiscrete limit of tau functions 
as in \cite{HHIKTT} for $A^{(1)}_{n-1}$ as a future problem.

\section{Tropical $R$ for $D^{(1)}_n$}\label{sec:TropicalR}
Let us recall the tropical $R$ for $D^{(1)}_n$ {}from \cite{KOTY}.
Let $x = (x_1,\ldots,x_n,\xb_{n-1},\ldots,\xb_1)$ be a set of variables.
The geometric $D^{(1)}_n$-crystal is the set $\B = \{x \}$ 
equipped with additional structures such as 
$\vphi_i, \veps_i$ and the geometric Kashiwara operators $e^c_i$.
The product $\B \times \cdots \times \B$ again becomes a 
geometric crystal, which is an analogue of the tensor product of 
crystals.
The tropical $R$ is a birational map 
$\B \times \B \rightarrow \B \times \B$ commuting with 
the geometric Kashiwara operators.
Leaving the precise definition to \cite{KOTY}, 
we concentrate here on its explicit form.

Set 
\begin{equation}\label{eq:level}
\ell(x) = x_1x_2\cdots x_n\xb_{n-1}\cdots \xb_2\xb_1,
\end{equation}
and call it the {\em level} of $x$.
On $\B$ we introduce the involutions 
$\sigma_1$ and $\sigma_n$ which preserve the level.
Explicitly for $x = (x_1,\ldots,x_n,\xb_{n-1},\ldots,\xb_1) \in \B$, 
they are defined by
\begin{align}
\label{eq:sigma1o}
\sigma_1 &: x_1 \longleftrightarrow \xb_1, \\
\label{eq:sigmano}
\sigma_n &: x_{n-1} \rightarrow x_{n-1}x_n, \quad
             \xb_{n-1} \rightarrow \xb_{n-1}x_n, \quad
             x_n \rightarrow 1/x_n,
\end{align}
where the variables not appearing in the above are unchanged.
On $\B \times \B$ we introduce $\sigma_1, \sigma_n$ analogously as
\begin{equation}\label{eq:sigma1n}
\sigma_a(x,y) = (\sigma_a(x),\sigma_a(y))\quad a=1,n.
\end{equation}
Furthermore we introduce the involution $\sigma_\ast$ by 
\begin{align}
\label{eq:ast}
\sigma_\ast &: x_i \longleftrightarrow \yb_i, \quad
\xb_i \longleftrightarrow y_i \quad (1 \leq i \leq n-1),
\quad x_n \longleftrightarrow y_n
\end{align}
for the elements $x=(x_1,\ldots,\xb_1)$ and
$y=(y_1,\ldots,\yb_1)$ of $\B$.
Note that $\sigma_1, \sigma_n$ and 
$\sigma_\ast$ are all commutative, among which 
$\sigma_\ast$ is the only one that 
mixes $x$ and $y$ and interchanges the levels.
Given a function $F = F(x,y)$, we will write 
$F^{\sigma_a} = F(\sigma_a(x,y))$ for 
$a=1,n, \ast$.

We set 
\begin{align}
\label{eq:16}
V_0 &= \ell (x) \frac{y_1}{\yb_1} + \ell (x) \sum_{m=2}^{n-1}
\left( \prod_{i=1}^{m-1} \frac{y_i}{x_i} \right)
\left( 1+\frac{y_{m}}{\yb_{m}} \right) +
x_n y_n \prod_{i=1}^{n-1} \xb_i y_i \\
\nonumber
&+ \ell (y) \frac{\xb_1}{x_1} + \ell (y) \sum_{m=2}^{n-1}
\left( \prod_{i=1}^{m-1} \frac{\xb_i}{\yb_i} \right)
\left( 1+\frac{\xb_{m}}{x_{m}} \right) +
\prod_{i=1}^{n-1} \xb_i y_i,
\end{align}
and define $V_i \, (1 \leq i \leq n-1)$ by the recursion relations:
\begin{equation}\label{eq:vind}
\begin{split}
&V_i = \yb_i\left(\frac{1}{\xb_i}V_{i-1} + (\ell (x) - \ell (y))
\Bigl(\frac{1}{x_i}+\frac{1}{\yb_i} \Bigr) \right),
\quad (1 \leq i \leq n-2)\\
&V_{n-1} = \frac{\yb_{n-1}}{y_n}\left(
\frac{y_n}{\xb_{n-1}}V_{n-2} +
(\ell (x) - \ell (y)) \Bigl(\frac{y_n}{x_{n-1}} +
\frac{1}{\yb_{n-1}}\Bigr)\right).
\end{split}
\end{equation}
In terms of $V_i$ we also define $U_i \, (1 \leq i \leq n-1)$ by
\begin{equation}\label{eq:w}
\begin{split}
& U_1 = V_0V^{\sigma_1}_0, \quad U_{n-1} = V_{n-1} V_{n-1}^{{\sigma_\ast}},
\\
& U_i = \left( \frac{1}{x_i} + \frac{1}{\yb_i} \right)^{-1}
\left(\frac{1}{y_i} V_i V_{i-1}^{\sigma_\ast} 
+ \frac{1}{\xb_i} V_{i-1} V_i^{\sigma_\ast}\right)
\quad (2 \leq i \leq n-2). 
\end{split}
\end{equation}
The tropical $R$ is a birational map $R(x,y) = (x',y')$ 
on $\B \times \B$ given by
\begin{equation}\label{eq:tropRform}
\begin{split}
& x_1' = y_1 \frac{V_0^{\sigma_1}}{V_1}, \quad
\xb_1' = \yb_1 \frac{V_0}{V_1}, \\
& x_i' = y_i \frac{V_{i-1}U_i}{V_iU_{i-1}}, \quad
\xb_i' = \yb_i \frac{V_{i-1}}{V_i}, \quad (2 \leq i \leq n-1)\\
& x_n' = y_n \frac{V_{n-1}}{V_{n-1}^{{\sigma_n}}},  \\
& y_1' = x_1 \frac{V_0}{V_1^{\sigma_\ast}}, \quad
\yb_1' = \xb_1 \frac{V_0^{\sigma_1}}{V_1^{\sigma_\ast}}, \\
& y_i' = x_i \frac{V_{i-1}^{\sigma_\ast}}{V_i^{\sigma_\ast}}, \quad
\yb_i' = \xb_i \frac{V_{i-1}^{\sigma_\ast} U_i}{V_i^{\sigma_\ast} U_{i-1}}, 
\quad (2 \leq i \leq n-1)\\
& y_n' = x_n \frac{V_{n-1}^{\sigma_n}}{V_{n-1}},
\end{split}
\end{equation}
In \cite{KOTY}, $U_i$ here was denoted by $W_i$ in Definition 4.9, 
and their transformation property under $\sigma_a$
was summarized in Table 1.
As shown therein, the tropical $R$ is subtraction-free (totally positive), 
interchanges the levels, commutes with the geometric Kashiwara operators, 
and satisfies the inversion relation 
$R(R(x,y)) = (x,y)$ and the Yang-Baxter relation.

\begin{remark}\label{rem:RAC}
Let $(x',y') = R(x,y)$ under the $D^{(1)}_n$ tropical $R$.
It is easily seen that 
$x_n=y_n=1$ is equivalent to $x'_n = y'_n=1$.
Similarly, $(x_1,y_1) = (\xb_1,\yb_1)$ is equivalent to 
$(x'_1,y'_1) = (\xb'_1,\yb'_1)$.
\end{remark}

\section{Bilinear equations}\label{sec:BLE}

Here we introduce a system of bilinear equations 
and study its properties such as 
existence and uniqueness of the solution.
It is a preparation for Section \ref{sec:B-tropR}, where 
the equations are related to our tropical $R$.

In the rest of the paper we fix the elements 
$\la = (\la_1,\ldots, \lab_1)$ and $\ka = (\ka_1,\ldots, \kab_1)$ of $\B$, 
and set 
\begin{equation}\label{eq:lkdef}
l = \ell(\la),\qquad k = \ell(\ka).
\end{equation}
We assume $lk(l-k) \neq 0$ throughout.
We introduce the variables which we call the {\em tau functions}:
\begin{equation*}
S_i, W_i, N_i, E_i\quad (1 \le i \le n-2),\qquad
\tau^0_i, \tau^1_i, \tau^2_i, \tau^3_i, \tau^4_i \quad (0 \le i \le n).
\end{equation*}
They will be called {\em generic} if they are nonzero.
Let $\alpha$ and $\beta$ be arbitrary nonzero parameters.
Our bilinear equations read as follows:
\begin{align*}
\kab_1N_1\tau^4_1 - \lab_1E_1\tau^2_1 &= \alpha\tau^0_0 \tau^1_0,\\
\ka_1N_1 \tau^4_0 - \la_1E_1\tau^2_0 &= \alpha\tau^0_1\tau^1_1,\\
\ka_iE_{i-1}N_i - \la_iN_{i-1}E_i &= \alpha\tau^0_i\tau^1_i\quad (2 \le i \le n-2),\\
\ka_{n-1}\ka_nE_{n-2}\tau^2_n - 
\la_{n-1}\la_nN_{n-2}\tau^4_n &= \alpha\tau^0_{n-1}\tau^1_n,\\
\ka_{n-1}E_{n-2}\tau^2_{n-1} - 
\la_{n-1}N_{n-2}\tau^4_{n-1} &= \alpha\tau^0_n\tau^1_{n-1}.
\end{align*}
The $i$-th equation here will be referred as 
$\langle 1,i\rangle\; (0 \le i \le n)$.
\begin{align*}
\ka_1N_1\tau^3_1 + \lab_1W_1\tau^1_1 &= \beta\tau^0_0\tau^2_0,\\
\kab_1N_1\tau^3_0 + \la_1W_1\tau^1_0 &= \beta\tau^0_1\tau^2_1,\\
\kab_iW_{i-1}N_i + \la_iN_{i-1}W_i &= \beta\tau^0_i\tau^2_i\quad (2 \le i \le n-2),\\
\ka_n\kab_{n-1}W_{n-2}\tau^1_{n-1} + \la_{n-1}N_{n-2}\tau^3_{n-1} 
&= \beta\tau^0_{n-1}\tau^2_{n-1},\\
\kab_{n-1}W_{n-2}\tau^1_n + \la_{n-1}\la_nN_{n-2}\tau^3_n 
&= \beta\tau^0_n\tau^2_n.
\end{align*}
The $i$-th equation here will be referred as 
$\langle 2,i\rangle\; (0 \le i \le n)$.
\begin{align*}
\ka_1S_1\tau^2_1 - \la_1W_1\tau^4_1 &= \alpha\tau^0_0 \tau^3_0,\\
\kab_1S_1 \tau^2_0 - \lab_1W_1\tau^4_0 &= \alpha\tau^0_1\tau^3_1,\\
\kab_iW_{i-1}S_i - \lab_iS_{i-1}W_i &= \alpha\tau^0_i\tau^3_i\quad (2 \le i \le n-2),\\
\kab_{n-1}W_{n-2}\tau^4_n - 
\lab_{n-1}S_{n-2}\tau^2_n &= \alpha\tau^0_{n-1}\tau^3_n,\\
\ka_n\kab_{n-1}W_{n-2}\tau^4_{n-1} - 
\la_n\lab_{n-1}S_{n-2}\tau^2_{n-1} &= \alpha\tau^0_n\tau^3_{n-1}.
\end{align*}
The $i$-th equation here will be referred as 
$\langle 3,i\rangle\; (0 \le i \le n)$.
\begin{align*}
\kab_1S_1\tau^1_1 + \la_1E_1\tau^3_1 &= \beta\tau^0_0\tau^4_0,\\
\ka_1S_1\tau^1_0 + \lab_1E_1\tau^3_0 &= \beta\tau^0_1\tau^4_1,\\
\ka_iE_{i-1}S_i + \lab_iS_{i-1}E_i &= \beta\tau^0_i\tau^4_i\quad (2 \le i \le n-2),\\
\ka_{n-1}E_{n-2}\tau^3_{n-1} + \la_n\lab_{n-1}S_{n-2}\tau^1_{n-1} 
&= \beta\tau^0_{n-1}\tau^4_{n-1},\\
\ka_{n-1}\ka_nE_{n-2}\tau^3_n + \lab_{n-1}S_{n-2}\tau^1_n 
&= \beta\tau^0_n\tau^4_n.
\end{align*}
The $i$-th equation here will be referred as 
$\langle 4,i\rangle\; (0 \le i \le n)$.

We call  $\{\langle J,i\rangle \mid 1 \le J \le 4, 0 \le i \le n\}$ the 
$D^{(1)}_n$ bilinear equations.
\begin{lemma}\label{lem:body}
For generic input $N_1,\ldots, N_{n-2}, W_1,\ldots, W_{n-2}$ 
and $\tau^J_0,\ldots, \tau^J_n$ with $J=1,2,3$, 
there exists a unique solution to the equations 
$\langle J,i\rangle (1 \le J \le 4, 1 \le i \le n-1)$ on the variables
$S_1,\ldots, S_{n-2}, \; E_1,\ldots, E_{n-2},\;
\tau^0_1,\ldots, \tau^0_{n-1}$ and 
$\tau^4_0,\ldots, \tau^4_n$.
\end{lemma}

\begin{proof}
The equations $\langle 2,1\rangle-\langle 2,n-1\rangle$ 
fix $\tau^0_1,\ldots, \tau^0_{n-1}$.
Regard $\langle 1,i\rangle$ and $\langle 3,i\rangle$ with $1 \le i \le n-1$ as 
linear equations on the unknowns 
$S_1,\ldots, S_{n-2}$, $E_1,\ldots, E_{n-2}$, $\tau^4_0$ and $\tau^4_n$.
Arrange them in a matrix form $A {\bf x} = \alpha{\bf b}$, where the 
column vectors ${\bf x}$ and ${\bf b}$ are specified by
\begin{align*}
&{\bf x}={}^t(\tau^4_0,E_1,\ldots,E_{n-2},\tau^4_n,S_{n-2},\ldots,S_1),\\
&{\bf b} = {}^t(\tau^0_1\tau^1_1,\ldots,\tau^0_{n-2}\tau^1_{n-2},
\tau^0_{n-1}\tau^1_n,\tau^0_{n-1}\tau^3_n,\tau^0_{n-2}\tau^3_{n-2},\ldots,
\tau^0_1\tau^3_1).
\end{align*}
Then one finds that the $2n-2$ by $2n-2$ coefficient matrix $A$ 
is an upper triangular Jacobi matrix except the bottom left element $-\lab_1W_1$.
Hence it is easy to see 
$\det A = (k-l)\tau^2_0\tau^2_n\prod_{i=1}^{n-2}(N_iW_i)$, which is 
nonzero for generic input.
Therefore the unknowns in ${\bf x}$ are uniquely determined.
With $\tau^0_1,\ldots, \tau^0_{n-1}$ and ${\bf x}$ at hand, one can determine 
the remaining ones $\tau^4_1,\ldots,\tau^4_{n-1}$ {}from 
$\langle 4,1\rangle-\langle 4,n-1\rangle$. 
\end{proof}

According to Lemma \ref{lem:body}, 
the bilinear equations $\langle J,0\rangle$ and 
$\langle J,n\rangle$ are not needed to determine the tau functions 
that will appear in the parameterization 
(\ref{eq:xy-tau})-(\ref{eq:xyprime-tau}).
However, they are essential 
in order to incorporate the involutions $\sigma_a$ 
in the previous section into tau functions.
The rest of this section concerns this point, leading to 
Propositions \ref{pr:unique-existence}, \ref{pr:invariance}, 
and ultimately to Proposition \ref{pr:involution}.

\begin{lemma}\label{lem:null12}
Suppose that $S_i,W_i,N_i,E_i\; (1 \le i \le n-2)$ and 
$\tau^J_i\;(0 \le J \le 4, 0 \le i \le n)$ 
satisfy $\langle J,i\rangle (1 \le J \le 4, 1 \le i \le n-1)$.
Then we have 
\begin{align}
&\alpha\lab_1\kab_1\tau^3_0\tau^1_1 + 
\beta\la_1\kab_1\tau^2_0\tau^4_1 -
\alpha\la_1\ka_1\tau^1_0\tau^3_1 - 
\beta\lab_1\ka_1\tau^4_0\tau^2_1 = 0, \label{eq:null-0}\\
&\alpha\tau^3_{n-1}\tau^1_n +
\beta\la_n\tau^4_n\tau^2_{n-1} - 
\alpha\la_n\ka_n\tau^1_{n-1}\tau^3_n -
\beta\ka_n\tau^2_n\tau^4_{n-1} = 0.\label{eq:null-n}
\end{align}
\end{lemma}

\begin{proof}
Write $\langle 1,1\rangle-\langle 4,1\rangle$ in the matrix form:
\begin{equation*}
\begin{pmatrix}
\ka_1\tau^4_0 & -\la_1\tau^2_0 & 0 & 0 \\
0 & \lab_1\tau^3_0 & \ka_1\tau^1_0 & 0 \\
0 & 0 & \kab_1\tau^2_0 & -\lab_1\tau^4_0 \\
\kab_1\tau^3_0 & 0 & 0 & \la_1\tau^1_0
\end{pmatrix}
\begin{pmatrix}N_1 \\ E_1 \\ S_1 \\ W_1 \end{pmatrix}
= \tau^0_1
\begin{pmatrix}\alpha\tau^1_1 \\ \beta\tau^4_1
\\ \alpha\tau^3_1 \\ \beta\tau^2_1 \end{pmatrix}.
\end{equation*}
The matrix on the left hand side is annihilated by
multiplying the row vector 
$(\lab_1\kab_1\tau^3_0, \la_1\kab_1\tau^2_0, 
-\la_1\ka_1\tau^1_0,-\lab_1\ka_1\tau^4_0)$ {}from the left.
So the same should happen on the right hand side, 
proving (\ref{eq:null-0}).
The relation (\ref{eq:null-n}) is verified similarly by using 
$\langle 1,n-1\rangle-\langle 4,n-1\rangle$.
\end{proof}

\begin{lemma}\label{lem:around-0}
Suppose that $S_i,W_i,N_i,E_i\; (1 \le i \le n-2)$, 
$\tau^J_i\;(1 \le J \le 4, 0 \le i \le n)$ and 
$\tau^0_i\;(1 \le i \le n-1)$ satisfy 
the equations $\langle J,i\rangle (1 \le J \le 4, 1 \le i \le n-1)$.
Then there exists a unique solution to the equations 
$\langle J,0\rangle\, (1 \le J \le 4)$ on $\tau^0_0$.
\end{lemma}

\begin{proof}
It suffices to check compatibility of 
the four equations $\langle J,0\rangle\, (1 \le J \le 4)$, that is, 
they all lead to the same $\tau^0_0$.
Multiply $\langle 1,0\rangle$ by $\tau^3_0$, and 
$\langle 3,0\rangle$ by $\tau^1_0$.
The difference of the resulting left hand sides can be grouped as
\begin{equation*}
\tau^4_1(\kab_1N_1\tau^3_0 + \la_1W_1\tau^1_0) - 
\tau^2_1(\lab_1E_1\tau^3_0 + \ka_1S_1\tau^1_0).
\end{equation*}
Due to $\langle 2,1\rangle$ and $\langle 4,1\rangle$ this vanishes, proving the 
compatibility of $\langle 1,0\rangle$ and $\langle 3,0\rangle$.
The compatibility of $\langle 2,0\rangle$ 
and $\langle 4,0\rangle$ is similarly shown 
by using $\langle 1,1\rangle$ and $\langle 3,1\rangle$. 
Let us check the compatibility of $\langle 2,0\rangle$ and $\langle 3,0\rangle$.
We are to show 
\begin{equation*}
\alpha\tau^3_0(\ka_1N_1\tau^3_1+\lab_1W_1\tau^1_1) - 
\beta\tau^2_0(\ka_1S_1\tau^2_1 - \la_1W_1\tau^4_1) = 0
\end{equation*}
Upon multiplying $\kab_1$, the first 
term contains the factor $\kab_1N_1\tau^3_0$ and the third 
does $\kab_1S_1\tau^2_0$. 
Eliminate them by using $\langle 2,1\rangle$ 
and $\langle 3,1\rangle$, respectively.
After a cancellation the resulting 
expression becomes $W_1$ times 
the left hand side of (\ref{eq:null-0}).
\end{proof}

\begin{lemma}\label{lem:around-n}
Suppose that $S_i,W_i,N_i,E_i\; (1 \le i \le n-2)$, 
$\tau^J_i\;(1 \le J \le 4, 0 \le i \le n)$ and 
$\tau^0_i\;(1 \le i \le n-1)$ satisfy 
the equations $\langle J,i\rangle (1 \le J \le 4, 1 \le i \le n-1)$.
Then there exists a unique solution to the equations 
$\langle J,n\rangle\, (1 \le J \le 4)$ on $\tau^0_n$.
\end{lemma}

\begin{proof}
This is verified similarly to Lemma \ref{lem:around-0} 
by using (\ref{eq:null-n}).
\end{proof}

Since the compatibility of the two systems 
$\langle 1,0\rangle-\langle 4,0\rangle$ and 
$\langle 1,n\rangle-\langle 4,n\rangle$ is trivial, 
Lemmas \ref{lem:body}, \ref{lem:around-0} and 
\ref{lem:around-n} lead to
\begin{proposition}[Unique existence]\label{pr:unique-existence}
For generic input $N_1,\ldots, N_{n-2}$, $W_1,\ldots, W_{n-2}$ 
and $\tau^J_0,\ldots, \tau^J_n$ with $J=1,2,3$, 
there exists a unique solution to 
the bilinear equations $\langle J,i\rangle (1 \le J \le 4, 0 \le i \le n)$ 
on the remaining variables
$S_1,\ldots, S_{n-2}$, $E_1,\ldots, E_{n-2}$ and
$\tau^J_0,\ldots, \tau^J_n\; (J=0,4)$.
\end{proposition}

Let us proceed to the automorphism of the 
bilinear equations.
For brevity the array 
${\mathcal T} = (\la,\ka ; S_i, W_i, N_i, E_i, \tau^J_i)$  
consisting of all the tau functions together with $(\la,\ka) \in \B\times \B$ 
will be called {\em data}.
We extend the involutions  $\sigma_a\, (a=1,n,\ast)$ on $\B \times \B$ 
defined in the previous section to the data as follows:
\begin{align}
&\sigma_a \; (a=1,n,\ast) \,\hbox{ acts on } \, (\la,\ka) \in \B\times \B \, 
\hbox{ by (\ref{eq:sigma1o})-(\ref{eq:ast}) as before},\label{eq:asbefore}\\
&\sigma_1: \; \tau^J_0 \longleftrightarrow \tau^J_1 \;\; (0 \le J \le 4),
\label{eq:sig1-tau}\\
&\sigma_n: \; \tau^J_{n-1} \longleftrightarrow \tau^J_n \;\; (0 \le J \le 4),
\label{eq:sign-tau}\\
&\; \sigma_\ast : \; 
\tau^0_{n-1} \longleftrightarrow \tau^0_n, \quad
\tau^2_{n-1} \longleftrightarrow \tau^2_n,\quad 
\tau^4_{n-1} \longleftrightarrow \tau^4_n, \label{eq:*-tau}\\
&\qquad \; \tau^1_i \longleftrightarrow \tau^3_i\;\; (0 \le i \le n-2),\quad
\tau^1_{n-1} \longleftrightarrow \tau^3_n,\quad
\tau^1_{n} \longleftrightarrow \tau^3_{n-1},\nonumber \\
&\qquad \; W_i \longleftrightarrow N_i,\quad \; 
S_i \longleftrightarrow E_i\;\; (1 \le i \le n-2).\nonumber
\end{align}
In each involution, the tau functions not appearing in the above 
are unchanged.
See also Proposition \ref{pr:involution}.
It is immediate to check 
\begin{proposition}[Invariance]\label{pr:invariance}
The involutions $\sigma_a\; (a=1,n,\ast)$ defined in 
(\ref{eq:asbefore})-(\ref{eq:*-tau}) on the data ${\mathcal T}$ 
are automorphism of 
the bilinear equations $\langle J,i\rangle (1 \le J \le 4, 0 \le i \le n)$.
\end{proposition}

In particular, the composition $\sigma_n\sigma_\ast$ 
interchanges $\langle 1,i \rangle$ and $\langle 3, i \rangle$
leaving $\langle 2,i \rangle$ and $\langle 4,i \rangle$ invariant.

\section{Bilinearization of Tropical $R$}\label{sec:B-tropR}

Given $\mu = (\mu_1,\ldots, \overline{\mu}_1) \in \B$ and the arrays 
\begin{equation*}
C = (C_1,\ldots,C_{n-2}), \quad \tau = (\tau_0,\ldots,\tau_n), 
\quad \tau' = (\tau'_0,\ldots,\tau'_n), 
\end{equation*}
we define
\begin{equation}\label{eq:def[]}
\begin{split}
[\mu; \tau, C, \tau'] &= (z_1,\ldots, z_n,\ol{z}_{n-1},\ldots,\ol{z}_1) \in \B,\\
z_1 &= \mu^{-1}_1\frac{\tau_0\tau'_1}{C_1},\quad 
\ol{z}_1 = \ol{\mu}^{-1}_1\frac{\tau_1\tau'_0}{C_1},\\
z_2 &= \mu^{-1}_2\frac{C_1\tau'_2}{C_2\tau'_0\tau'_1},\quad 
\ol{z}_2 = \ol{\mu}^{-1}_2\frac{C_1\tau_2}{C_2\tau_0\tau_1},\\
z_i &= \mu^{-1}_i\frac{C_{i-1}\tau'_i}{C_i\tau'_{i-1}},\quad
\ol{z}_i = \ol{\mu}^{-1}_i\frac{C_{i-1}\tau_i}{C_i\tau_{i-1}}
\quad (3 \le i \le n-2),\\
z_{n-1} &= 
\mu^{-1}_{n-1}\frac{C_{n-2}\tau'_{n-1}}{\tau_{n-1}\tau'_{n-2}},\quad
\ol{z}_{n-1} = \ol{\mu}^{-1}_{n-1}
\frac{C_{n-2}\tau_{n}}{\tau_{n-2}\tau'_{n}},\\
z_n &= \mu^{-1}_n\frac{\tau_{n-1}\tau'_n}{\tau_n\tau'_{n-1}}
\end{split}
\end{equation}
in case  $n \ge 4$. 
For $n=3$ we slightly modify it into
\begin{align*}
z_1 &= \mu^{-1}_1\frac{\tau_0\tau'_1}{C_1},\quad 
\ol{z}_1 = \ol{\mu}^{-1}_1\frac{\tau_1\tau'_0}{C_1},\\
z_2 &= \mu^{-1}_2\frac{C_1\tau'_2}{\tau_2\tau'_0\tau'_1},\quad 
\ol{z}_2 = \ol{\mu}^{-1}_2\frac{C_1\tau_3}{\tau_0\tau_1\tau'_3},\\
z_3 &= \mu^{-1}_3\frac{\tau_{2}\tau'_3}{\tau_3\tau'_{2}}.
\end{align*}
Note that $\ell([\mu;\tau,C,\tau']) = \ell(\mu)^{-1}$.
{}From the data ${\mathcal T}= (\la,\ka ; S_i, W_i, N_i, E_i, \tau^J_i)$,
we construct $x,y,x',y' \in \B$ as 
\begin{align}
&x = [\la; \tau^3, W, \tau^2],\quad 
y = [\ka; \tau^2, N, \tau^1],\label{eq:xy-tau}\\
&x' = [\ka; \tau^3, S, \tau^4],\quad
y' = [\la; \tau^4, E, \tau^1].\label{eq:xyprime-tau}
\end{align}
Note {}from (\ref{eq:lkdef}) that 
\begin{equation}\label{eq:xylevel}
\ell(x) = \ell(y') = l^{-1},\quad \ell(x') = \ell(y) = k^{-1}.
\end{equation}
Given any pairs $(x,y)$ and $(\la,\ka)$ such that 
$\ell(x) = \ell(\la)^{-1}$ and 
$\ell(y) = \ell(\ka)^{-1}$, 
the parameterization (\ref{eq:xy-tau}) is always possible due to 
\begin{proposition}\label{pr:adequate}
For any $z, \mu \in \B$ such that $\ell(z) = \ell(\mu)^{-1}$, 
the equation $[\mu; \tau,C,\tau'] = z$ with fixed $\tau$, (resp. $\tau'$)
on the variables $(C,\tau')$ (resp. $(C,\tau)$) 
admits a one parameter family of solutions.
\end{proposition}

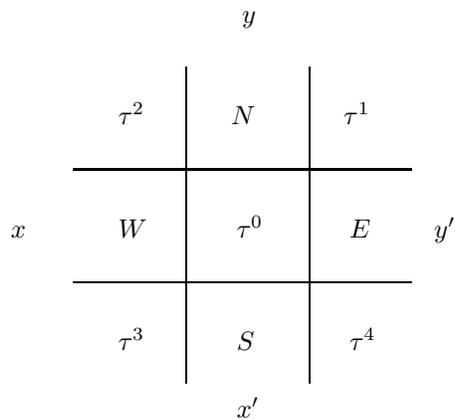
\begin{figure}[h]
\caption{Tau functions}\label{fig:fig1}
\setlength{\unitlength}{1.5mm}
\begin{picture}(45,40)(-9.5,-2)
\put(3,24){$\tau^{2}$}
\put(3,14){$W$}
\put(3,4){$\tau^{3}$}
\put(13,24){$N$}
\put(13.5,14){$\tau^{0}$}
\put(13.5,4){$S$}
\put(23,24){$\tau^{1}$}
\put(23.5,14){$E$}
\put(23.5,4){$\tau^{4}$}
\put(-1,20){\line(1,0){30}}
\put(29,10){\line(-1,0){30}}
\put(20,29){\line(0,-1){28}}
\put(9,1){\line(0,1){28}}
\put(-6.5,14){$x$}
\put(14,33){$y$}
\put(31,14){$y'$}
\put(13.5,-2){$x'$}
\end{picture}
\end{figure}

It is customary to attach a vertex diagram (cross) to the relation like
$R(x,y) = (x',y')$ for a quantum or combinatorial $R$, where 
the four edges correspond to $x,y,x',y'$.
Figure \ref{fig:fig1} is an analogue of such a diagram, which 
may be of help to recognize the pattern 
(\ref{eq:xy-tau})-(\ref{eq:xyprime-tau}).
Elements of $\B$ are represented there with double lines.
Consequently the vertex diagram is separated into several domains where 
various tau functions live.

It is readily seen that the
transformation ${\mathcal T} \rightarrow \sigma_a({\mathcal T})$ of the data
induces the change $(x,y) \rightarrow \sigma_a(x,y)$ and 
$(x',y') \rightarrow \sigma_a(x',y')$.
More generally we have 
\begin{proposition}\label{pr:involution}
Let $F = F(x,y)$ be any function on $\B \times \B$ and \\
${\mathcal T} = (\la,\ka ; S_i, W_i, N_i, E_i, \tau^J_i)$ be 
any data obeying all the bilinear equations.
Denote by $F\vert_{\mathcal T}$ an expression of $F$ 
resulting {}from the substitution of (\ref{eq:xy-tau}) 
followed by any possible application of the bilinear equations.
Then one has 
$F^{\sigma_a}\vert_{\mathcal T} = 
\sigma_a\left(F\vert_{\mathcal T}\right)\, (a=1,n,\ast)$, 
where $\sigma_a$ on the left hand side is 
specified by (\ref{eq:sigma1o})-(\ref{eq:ast}), whereas on the 
right hand side by 
(\ref{eq:asbefore})-(\ref{eq:*-tau}).
The same fact is valid also for (\ref{eq:xyprime-tau}) and 
$F = F(x',y')$.
\end{proposition}

\begin{proof}
This is due to Propositions \ref{pr:unique-existence} and 
\ref{pr:invariance}
\end{proof}

\begin{lemma}\label{lem:UV}
Suppose that the data $(\la,\ka ; S_i, W_i, N_i, E_i, \tau^J_i)$ satisfies
all the bilinear equations. Then under the substitution
(\ref{eq:xy-tau}), the functions $U_i, V_i$ 
defined in (\ref{eq:16})-(\ref{eq:w}) are expressed as follows:
\begin{align*}
&V_0 = V^{\sigma_\ast}_0 = \frac{\gamma\tau^2_0\tau^4_0}{\tau^1_0\tau^3_0}, 
\;\;
V^{\sigma_1}_0 = \frac{\gamma\tau^2_1\tau^4_1}{\tau^1_1\tau^3_1}, 
\;\;
U_1 = \frac{\gamma^2\tau^2_0\tau^2_1\tau^4_0\tau^4_1}
{\tau^1_0\tau^1_1\tau^3_0\tau^3_1},\\
&V_1 = V^{\sigma_1}_1 = 
\frac{\gamma S_1\tau^2_0\tau^2_1}{N_1\tau^3_0\tau^3_1}, 
\;\;
V^{\sigma_\ast}_1 = 
\frac{\gamma E_1\tau^2_0\tau^2_1}{W_1\tau^1_0\tau^1_1}, 
\\
&V_i = V^{\sigma_1}_i = \frac{\gamma S_i\tau^2_i}{N_i\tau^3_i}, 
\;\;
V^{\sigma_\ast}_i = \frac{\gamma E_i\tau^2_i}{W_i\tau^1_i}, 
\;\;
U_i =  \frac{\gamma^2\tau^2_i\tau^4_i}{\tau^1_i\tau^3_i}
\;\;(2 \le i \le n-2),\\
&V_{n-1} = V^{\sigma_1}_{n-1} = 
\frac{\gamma\tau^2_n\tau^4_n}{\tau^1_n\tau^3_n}, 
\;\;
V^{\sigma_\ast}_{n-1} = V^{\sigma_n}_{n-1} 
= \frac{\gamma\tau^2_{n-1}\tau^4_{n-1}}
{\tau^1_{n-1}\tau^3_{n-1}}, 
\;\; 
U_{n-1} = \frac{\gamma^2\tau^2_{n-1}\tau^2_{n}\tau^4_{n-1}\tau^4_{n}}
{\tau^1_{n-1}\tau^1_{n}\tau^3_{n-1}\tau^3_{n}},
\end{align*}
where $\gamma = (k-l)\beta/(lk\alpha)$.
\end{lemma}
The proof is available in Appendix \ref{sec:appA}.
The main result of this section is
\begin{theorem}[Bilinearization of tropical $R$]\label{th:bilinear}
Suppose that the data $(\la,\ka ; S_i, W_i, N_i$, $E_i, \tau^J_i)$ satisfies
all the bilinear equations.  Then $x,y,x',y'$ specified in 
(\ref{eq:xy-tau}) and (\ref{eq:xyprime-tau}) obey the relation
$R(x,y) = (x',y')$.
\end{theorem}
\begin{proof}
Substitute (\ref{eq:xy-tau}) into the right hand sides of (\ref{eq:tropRform}) 
and apply Lemma \ref{lem:UV}.
Then the result agrees with (\ref{eq:xyprime-tau}).
\end{proof}

Given $(x,y)$, the tau functions $\tau^1, \tau^2, \tau^3, N, W$ 
satisfying (\ref{eq:xy-tau}) are not unique.
Theorem \ref{th:bilinear} guarantees that 
$(x',y')$  in (\ref{eq:xyprime-tau}) is independent of their choice
under the bilinear equations.
{}From Proposition \ref{pr:involution} and the remark preceding it,
we find that $R(\sigma_a(x,y)) = \sigma_a(R(x,y))$ in agreement with 
\cite{KOTY}. 
We conclude that the tropical $R$ acts on the data as the interchange
\begin{equation*}
\la \longleftrightarrow \ka, \;\; 
\tau^2_i \longleftrightarrow \tau^4_i,\;\;
W_i \longleftrightarrow S_i,\;\;
N_i \longleftrightarrow E_i
\end{equation*}
for all $i$.
Combined with 
$\alpha \rightarrow -\alpha$, 
this is another automorphism of our bilinear equations
that interchanges $\langle 2,i \rangle$ and $\langle 4, i \rangle$
leaving $\langle 1,i \rangle$ and $\langle 3,i \rangle$ invariant.
This is a complementary transformation with the $\sigma_n\sigma_\ast$
mentioned after Proposition \ref{pr:invariance}.

\section{Solutions of bilinear equations}\label{sec:sol}

There is a link between our bilinear equations and 
the soliton theory in terms of tau functions and 
infinite dimensional Lie algebras \cite{JM}.
We show in Section \ref{subsec:local} 
that certain vacuum expectation values of the 
2 component free fermions provide solutions 
to the bilinear equations.
They are parametrized with the elements of the algebras 
$D'_\infty \simeq D_\infty$ and their reduction.
Then in Section \ref{subsec:tvertex} 
the solutions for the local equations 
are naturally extended to the tropical vertex model on 
the two dimensional square lattice where the
tropical $R$ plays the role of local time evolution 
at each vertex.
As for notations and basic facts on the free fermion approach, see
Appendix \ref{sec:appB}.

\subsection{Tau functions as vacuum expectation values}\label{subsec:local}

In this subsection, the letters $x,y$, etc. stand for the 
array of infinitely many time variables as in Appendix \ref{sec:appB}.
They should not be 
confused with elements in ${\mathcal B}$, which have now been  
effectively replaced by the tau functions.
Thus we refresh Figure \ref{fig:fig1} into the following:

\begin{figure}[h]
\caption{Increment time arrays}\label{fig:fig2}
\setlength{\unitlength}{1.5mm}
\begin{picture}(45,40)(-9.5,-2)
\put(-9,9.5){${\tilde \veps}(L^{-1})$}
\put(-9,19.5){$\veps(L^{-1})$}
\put(5,30.5){${\tilde \veps}(K^{-1})$}
\put(16,30.5){$\veps(K^{-1})$}
\put(3,24){$\tau^{2}$}
\put(3,14){$W$}
\put(3,4){$\tau^{3}$}
\put(13,24){$N$}
\put(13.5,14){$\tau^{0}$}
\put(13.5,4){$S$}
\put(23,24){$\tau^{1}$}
\put(23.5,14){$E$}
\put(23.5,4){$\tau^{4}$}
\put(-1,20){\line(1,0){30}}
\put(29,10){\line(-1,0){30}}
\put(20,29){\line(0,-1){28}}
\put(9,1){\line(0,1){28}}
\end{picture}
\end{figure}
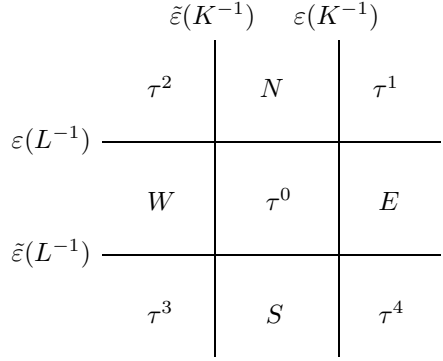

Here we have introduced the parameters $K,L$, 
which will play a role of level in ${\mathcal B}$.
Each line is assigned with the array $\veps$ or ${\tilde \veps}$ 
determined {}from these parameters according to 
(\ref{eq:eps}).
We assign the time variables $x^J$ 
on the 9 domains 
$J=0,1,2,3,4,N,S,W,E$ in Figure \ref{fig:fig2} so that the 
following recursion relations are satisfied:
\begin{equation}\label{eq:times}
\begin{split}
&x^W - x^3 = x^0 - x^S = x^E - x^4 = {\tilde \veps}(L^{-1}), \\
&x^2 - x^W = x^N - x^0 = x^1 - x^E = \veps(L^{-1}),\\
&x^N - x^2 = x^0-x^W = x^S - x^3 = {\tilde \veps}(K^{-1}),\\
&x^1 - x^N = x^E - x^0 = x^4 - x^S = \veps(K^{-1}).
\end{split}
\end{equation}
These relations determine all the $x^J$'s consistently 
upon specifying any one of them as an initial condition.
We fix $x^1$ to an arbitrary {\em odd} element 
in the sense of Appendix \ref{subsec:single}.
Then {}from the construction (\ref{eq:times}), 
all the time variables on the corner domains are odd, i.e., 
\begin{align}\label{eq:allodd}
{\widetilde {x^J}} = x^J\quad J = 1,2,3,4.
\end{align}
On the other hand the variables in the domains $N,W,S,E$ satisfy
\begin{equation}\label{eq:xshift}
\begin{split}
&{\widetilde {x^N}} = x^N + \veps(K^{-1}) - {\tilde \veps}(K^{-1}),
\quad 
{\widetilde {x^S}} = x^S + \veps(K^{-1}) - {\tilde \veps}(K^{-1}),\\
&{\widetilde {x^W}} = x^W + \veps(L^{-1}) - {\tilde \veps}(L^{-1}),
\quad 
{\widetilde {x^E}} = x^E + \veps(L^{-1}) - {\tilde \veps}(L^{-1}).
\end{split}
\end{equation}

For $i \ge 1$ define further 
\begin{align}
&x^J_i = x^J + z_i = (x^J_{i,1}, x^J_{i,2},x^J_{i,3},\ldots),
\nonumber\\
&z_1 = 0,\quad  z_i = -\sum_{k=2}^i\veps(a^{-1}_k) \;\; (i \ge 2),
\label{eq:zdef}
\end{align}
where $a_2, a_3, \ldots$ are nonzero parameters. 
Note that $x^J_1 = x^J$ and $x^J_i + \veps(a^{-1}_i) = x^J_{i-1}$.

Let $g' \in e^{D'_\infty}$ and 
$g \in e^{D_\infty}$ be the elements 
that are related as in (\ref{eq:ggp}). 
We shall use the vacuum expectation values 
(\ref{eq:fF}) exclusively for odd $y = (y_1,0,y_3,0,\ldots)$ with 
fixed $y_1, y_3,\ldots$.
Thus we simply write them as
$F_{l_1,l_2;l}(x;g)$ and $f_l(x;g')$.
Our solution to the 
bilinear equations is constructed in two steps.
In {\em Step 1}, we construct solutions in the 
infinite $n$ limit.
In {\em Step 2}, we will impose certain reduction condition on the 
elements $g, g'$ to make the equations close for finite $n$. 

{\em Step 1}.
We identify the tau functions with the vacuum expectation 
values as
\begin{align}
&X_j = F_{1,1}(x^X_j;g)  \;\;  j \ge 1 \quad (X=S,W,N,E),\label{eq:X}\\
&\tau_j^J = \begin{cases}
f_j(x^J_1;g') & j = 0,1, \\
F_{1,1}(x^J_j;g) & j \ge 2,
\end{cases}\quad (J = 1,2,3,4),\label{eq:tau1234}\\
&\tau^0_j = \begin{cases}
(F_{0,1} + (-1)^jiF_{0,1;1})(x^0_1;g) & j = 0,1, \\
F_{1,1}(x^0_{j-1};g) & j \ge 2.
\end{cases}\label{eq:tau0}
\end{align}
In the bilinear equations $\langle J,i \rangle$
in Section \ref{sec:BLE}, consider the $n$ infinity limit,
where one just forgets $\langle J,n-1 \rangle$ and 
$\langle J,n \rangle$ for all $J=1,2,3,4$.
We call the resulting system the $D^{(1)}_\infty$ 
bilinear equations.

\begin{lemma}\label{lem:0chikaku}
\begin{align*}
&F_{1,1}(x^J_1) = \tau^J_0\tau^J_1\quad J = 1,2,3,4,\\
&\left.\begin{array}{r}
2F_{0,1}(x^N_1;g) \\ 
2iF_{0,1;1}(x^N_1;g) 
\end{array}\right\}= 
\tau^1_1\tau^2_1\pm\tau^1_0\tau^2_0,\\
&\left.\begin{array}{r}
2F_{0,1}(x^W_1;g) \\ 
2iF_{0,1;1}(x^W_1;g) 
\end{array}\right\}= 
\tau^2_1\tau^3_1\pm\tau^2_0\tau^3_0,\\
&\left.\begin{array}{r}
2F_{0,1}(x^S_1;g) \\ 
2iF_{0,1;1}(x^S_1;g) 
\end{array}\right\}= 
\tau^3_1\tau^4_1\pm\tau^3_0\tau^4_0,\\
&\left.\begin{array}{r}
2F_{0,1}(x^E_1;g) \\ 
2iF_{0,1;1}(x^E_1;g) 
\end{array}\right\}= 
\tau^4_1\tau^1_1\pm\tau^4_0\tau^1_0.
\end{align*}
\end{lemma}
\begin{proof}
The first equality follows {}from Lemma \ref{lem:jm7.6} and 
(\ref{eq:allodd}).
The other relations are due to Lemma \ref{lem:sapix}
and (\ref{eq:xshift}).
\end{proof}

\begin{proposition}\label{pr:n-infinite}
The parameterization (\ref{eq:X})-(\ref{eq:tau0}) 
solves the $D^{(1)}_\infty$ bilinear equations 
with $(a_1=0)$
\begin{equation*}
\begin{split}
&\la_i = L - a_i, \; \lab_i = L+a_i\quad (i \ge 1),\\
&\ka_i = K - a_i, \; \kab_i = K+a_i\quad (i \ge 1),\\
&\alpha = K-L,\; \beta = K+L.
\end{split}
\end{equation*}
\end{proposition}

\begin{proof}
First consider $\langle J, i \rangle$ with $i \ge 2$.
Setting 
$l_1=l_2=1, l=0, x = x^0_i, b_3 = a_i$ in (\ref{eq:bl4}), we have
\begin{equation*}
\begin{split}
&(b_2-a_i)F_{1,1}(x^0_{i-1}+\veps(b^{-1}_2))
F_{1,1}(x^0_i+\veps(b^{-1}_1)) \\
&-(b_1-a_i)F_{1,1}(x^0_{i-1}+\veps(b^{-1}_1))
F_{1,1}(x^0_{i}+\veps(b^{-1}_2)) \\
&= (b_2-b_1)F_{1,1}(x^0_{i-1})
F_{1,1}(x^0_i+\veps(b^{-1}_1)+\veps(b^{-1}_2)),
\end{split}
\end{equation*}
where the dependence on $g \in e^{D_\infty}$ is suppressed.
Due to (\ref{eq:times}) 
this yields $\langle J,i \rangle$ with $1 \le J \le 4$ 
upon taking $(b_1,b_2) = (L,K), (L,-K)$, $(-L,-K)$, $(-L,K)$.
Next we treat $\langle J, 0 \rangle$ and 
$\langle J, 1 \rangle$.
Setting $l_1=l_2=1, l=0, x=x^0_1, b^{-1}_3=0$ in 
(\ref{eq:bl3}) $\mp i$ (\ref{eq:bl2}),
we have
\begin{align*}
&b_2F_{1,1}(x^0_1+\veps(b^{-1}_1))
(F_{0,1}\mp iF_{0,1;1})(x^0_1+\veps(b^{-1}_2))\\
&-b_1F_{1,1}(x^0_1+\veps(b^{-1}_2))
(F_{0,1}\mp iF_{0,1;1})(x^0_1+\veps(b^{-1}_1))\\
&= (b_2-b_1)F_{1,1}(x^0_1+\veps(b^{-1}_1)+\veps(b^{-1}_2))
(F_{0,1}\mp iF_{0,1;1})(x^0_1).
\end{align*}
By setting $(b_1,b_2) = (L,K)$, $(L,-K)$, $(-L,-K)$, $(-L,K)$, and 
applying Lemma \ref{lem:0chikaku}, this becomes 
$\langle J,1\rangle \tau^J_0$ and 
$\langle J,0\rangle \tau^J_1$ for $1 \le J \le 4$.
\end{proof}

{\em Step 2}.
We choose the elements $g' \in e^{D'_\infty}$ and 
$g \in e^{D_\infty}$ related by (\ref{eq:ggp}) as
\begin{align}
&g' = \exp\Bigl(\sum_{j=1}^N b_j\phi(p_j)\phi(q_j) 
+ \sum_{j=1}^M c_j\phi(p'_j){\hat \phi}(q'_j)\Bigr),\label{eq:gp}\\
&g = \exp\Bigl(\sum_{j=1}^N b_j\bigl(
\psi^{(1)}(p_j)\psi^{(1)\ast}(-q_j) -
\psi^{(1)}(q_j)\psi^{(1)\ast}(-p_j)\bigr) \label{eq:g}\\
&\qquad + \sum_{j=1}^M c_j\bigl(
\psi^{(1)}(p'_j)\psi^{(2)\ast}(-q'_j) - 
\psi^{(2)}(q'_j)\psi^{(1)\ast}(-p'_j)\bigr)\Bigr).\nonumber
\end{align}
On these elements we impose the reduction condition:
\begin{equation}\label{eq:red}
p_j^2A(p_j)A(-p_j) = q_j^2A(q_j)A(-q_j),\qquad 
{p'}^2_jA(p'_j)A(-p'_j) = {q'}^2_j
\end{equation}
for all $j$, where the function $A$ is defined by
\begin{equation*}
A(p) = \prod_{k=2}^{n-1}(1-\frac{p}{a_k})
\end{equation*}
in terms of $a_2, a_3, \ldots$ introduced in (\ref{eq:zdef}).

\begin{lemma}\label{lem:red}
Under the condition (\ref{eq:red}), 
there exists $h' \in e^{D'_\infty}$ and 
$h \in e^{D_\infty}$ that are related by $h = \iota(h')\kappa\iota(h')$ 
as in (\ref{eq:ggp}) and satisfy
\begin{equation*}
F_{l_1,l_2;l}(x+z_{n-1},y;\omega(g)) 
= F_{l_1,l_2;l}(x,y;h)
\end{equation*}
for any $x$ and odd $y$.
\end{lemma}
\begin{proof}
Under the map $\omega$ specified in (\ref{eq:omega}), 
$\psi^{(\alpha)}(p)$ and $\psi^{(\alpha)\ast}(q)$ 
are transformed into 
$p\psi^{(\alpha)}(p)$ and 
$q^{-1}\psi^{(\alpha)\ast}(q)$.
Then under the time evolution 
$\hbox{Ad}e^{H(z_{n-1},0)}$, they are further changed into
$pA(p)^{\delta_{\alpha 1}}\psi^{(\alpha)}(p)$ and 
$q^{-1}A(q)^{-\delta_{\alpha 1}}\psi^{(\alpha)\ast}(q)$, respectively.
Therefore if $g$ in (\ref{eq:g}) is denoted by $e^X$, we get 
$\hbox{Ad}e^{H(z_{n-1},0)}(\omega(g)) = e^{X'}$ with 
\begin{equation}\label{eq:xp}
\begin{split}
X' = &\sum_jb_j\Bigl(
-\frac{p_jA(p_j)}{q_jA(-q_j)}
\psi^{(1)}(p_j)\psi^{(1)\ast}(-q_j) +
\frac{q_jA(q_j)}{p_jA(-p_j)}
\psi^{(1)}(q_j)\psi^{(1)\ast}(-p_j)\Bigr) \\
&+ \sum_j c_j\Bigl(-\frac{p'_jA(p'_j)}{q'_j}
\psi^{(1)}(p'_j)\psi^{(2)\ast}(-q'_j) +
\frac{q'_j}{p'_jA(-p'_j)}
\psi^{(2)}(q'_j)\psi^{(1)\ast}(-p'_j)\Bigr).
\end{split}
\end{equation}
Under the condition (\ref{eq:red}), we find 
\begin{align*}
&X' = X\vert_{b_j \rightarrow b'_j, c \rightarrow c'_j},\\
&b'_j = -\frac{q_jA(q_j)}{p_jA(-p_j)}b_j,\quad 
c'_j = -\frac{p'_jA(p'_j)}{q'_j}c_j.
\end{align*}
In view of the definition (\ref{eq:fF}) and 
$H(x+z_{n-1},y) = H(x,y) + H(z_{n-1},0)$,
the elements $h'$ and $h$ in question are 
obtained by replacing $b_j, c_j$ by $b'_j, c'_j$ 
in $g'$ and $g$ given in (\ref{eq:gp}) and (\ref{eq:g}), respectively.
\end{proof}
Using $h'$ and $h$ in Lemma \ref{lem:red}, 
we now modify (\ref{eq:X})--(\ref{eq:tau0}) 
in {\em Step 1} into the finite $n$ version as
\begin{align}
&X_j = F_{1,1}(x^X_j;g)  \;\;  1 \le j \le n-2 
\quad (X=S,W,N,E),\label{eq:mX}\\
&\tau_j^J = \begin{cases}
f_j(x^J_1;g') & j = 0,1, \\
F_{1,1}(x^J_j;g) & 2 \le j \le n-2, \\
f_{n-j}(x^J_1;h') & j = n-1,n,
\end{cases}\quad (J = 1,2,3,4),\label{eq:mtau1234}\\
&\tau^0_j = \begin{cases}
(F_{0,1} + (-1)^jiF_{0,1;1})(x^0_1;g) & j = 0,1, \\
F_{1,1}(x^0_{j-1};g) & 2 \le j \le n-2,\\
(F_{0,0} + 
\frac{(-1)^{n-j}ia_{n-1}}{KL}F_{-1,1;1})
(x^0_{n-2};\omega(g))& j = n-1,n.\\
\end{cases}\label{eq:mtau0}
\end{align}
\begin{lemma}\label{lem:nchikaku}
\begin{align*}
&F_{0,0}(x^J_{n-1};\omega(g)) = \tau^J_{n-1}\tau^J_n\quad 
J = 1,2,3,4,\\
&\left.\begin{array}{r}
2F_{0,0}(x^N_{n-1};\omega(g))\\ 
2iK^{-1}F_{-1,1;1}(x^N_{n-1};\omega(g)) 
\end{array}\right\}=
\tau^1_n\tau^2_{n-1}\pm\tau^2_n\tau^1_{n-1},\\
&\left.\begin{array}{r}
2F_{0,0}(x^W_{n-1};\omega(g))\\ 
2iL^{-1}F_{-1,1;1}(x^W_{n-1};\omega(g)) 
\end{array}\right\}=
\tau^2_n\tau^3_{n-1}\pm\tau^3_n\tau^2_{n-1},\\
&\left.\begin{array}{r}
2F_{0,0}(x^S_{n-1};\omega(g))\\ 
2iK^{-1}F_{-1,1;1}(x^S_{n-1};\omega(g)) 
\end{array}\right\}=
\tau^4_n\tau^3_{n-1}\pm\tau^3_n\tau^4_{n-1},\\
&\left.\begin{array}{r}
2F_{0,0}(x^E_{n-1};\omega(g))\\ 
2iL^{-1}F_{-1,1;1}(x^E_{n-1};\omega(g)) 
\end{array}\right\}=
\tau^1_n\tau^4_{n-1}\pm\tau^4_n\tau^1_{n-1}.
\end{align*}
\end{lemma}
\begin{proof}
The first relation follows {}from Lemma \ref{lem:red}, 
Lemma \ref{lem:jm7.6} and (\ref{eq:allodd}).
{}From Lemma \ref{lem:red}, 
the left hand side of the second equation 
$2F_{0,0}(x^N_{n-1};\omega(g))$ is equal to $2F_{0,0}(x^N_1;h)$.
{}From (\ref{eq:times}), (\ref{eq:xshift}) and 
Lemma \ref{lem:sapix}, we get
\[
2F_{0,0}(x^N_1;h) = f_0(x^1_1;h')f_1(x^2_1;h')+
f_0(x^2_1;h')f_1(x^1_1;h'),
\]
which is the right hand side.
The other relations can be checked similarly.
\end{proof}
\begin{theorem}\label{th:n-finite}
The parameterization (\ref{eq:mX})-(\ref{eq:mtau0}) 
solves the $D^{(1)}_n$ bilinear equations 
with $(a_1=0)$
\begin{equation*}
\begin{split}
&\la_i = L - a_i, \; \lab_i = L+a_i
\quad (1 \le i \le n-1),\; \la_n=1,\\
&\ka_i = K - a_i, \; \kab_i = K+a_i
\quad (1 \le i \le n-1),\; \ka_n=1,\\
&\alpha = K-L,\; \beta = K+L.
\end{split}
\end{equation*}
\end{theorem}
\begin{proof}
Thanks to Proposition \ref{pr:n-infinite} we are left to show 
$\langle J, n-1\rangle$ and $\langle J, n\rangle$ only.
Let us illustrate the proof for $J=1$ since  
the equations for the other $J$ follow {}from the same argument by 
taking the parameters $b_1, b_2$ in the proof similarly to 
Proposition \ref{pr:n-infinite}.
In (\ref{eq:bl4}) set $l_1=l_2=l=0$, $(b_1,b_2) = (L,K)$,
$b_3=a_{n-1}$, $x=x^0_{n-1}$ and $g \rightarrow \omega(g)$.
The result reads
\begin{equation*}
\begin{split}
&\ka_{n-1}F_{0,0}(x^E_{n-2};\omega(g))
F_{0,0}(x^N_{n-1};\omega(g)) - 
\la_{n-1}F_{0,0}(x^N_{n-2};\omega(g))
F_{0,0}(x^E_{n-1};\omega(g))\\
&=(K-L)F_{0,0}(x^1_{n-1};\omega(g))
F_{0,0}(x^0_{n-2};\omega(g)).
\end{split}
\end{equation*}
Owing to (\ref{eq:sym2}), we know 
$F_{0,0}(x^E_{n-2};\omega(g)) = F_{1,1}(x^E_{n-2};g) = E_{n-2}$ 
and similarly 
$F_{0,0}(x^N_{n-2};\omega(g)) = N_{n-2}$.
Rewriting $F_{0,0}(x^1_{n-1};\omega(g)), F_{0,0}(x^N_{n-1};\omega(g))$ and 
$F_{0,0}(x^E_{n-1};\omega(g))$ by Lemma \ref{lem:nchikaku}, 
we find that the above equation coincides with the combination
$\langle 1, n-1\rangle\tau^1_{n-1} + 
\langle 1, n\rangle\tau^1_n$.
The other combination 
$\langle 1, n-1\rangle\tau^1_{n-1} - 
\langle 1, n\rangle\tau^1_n$ can be shown 
by using (\ref{eq:bl1}) with $l_1=l_2=l=0$ similarly.
\end{proof}

\subsection{Tropical vertex model}\label{subsec:tvertex}
Consider the two dimensional square lattice ${\mathcal L'}$
equipped with the space-time coordinates $(s,t) \in \Z^2$.
The coordinate $s$ (resp. $t$) 
increases rightward (resp. downward) and each vertical 
(resp. horizontal) line corresponds to
$s=$constant (resp. $t=$constant).
On each edge of ${\mathcal L'}$ 
we assign an element in ${\mathcal B}$ so that 
those four surrounding a vertex obey 
the relation $R(x,y) = (x',y')$, where 
$x,y,x',y'$ are the ones on the east, north, south and 
west edges.
We call the resulting 
two dimensional system the {\em tropical vertex model}.
Unlike the usual vertex models in statistical mechanics \cite{B},
it is a deterministic system in the sense that 
all the edge variables are determined uniquely {}from their 
values on the northwest or southeast boundaries of ${\mathcal L'}$.
The tropical vertex model reduces to the 
$D^{(1)}_n$-soliton cellular automaton 
\cite{HKT1,HKT2,HKOTY} in the ultradiscrete limit, where 
the tropical $R$ is replaced by the combinatorial $R$ \cite{HKOT} and 
${\mathcal B}$ by the $D^{(1)}_n$ crystal \cite{KKM}.
{}From (\ref{eq:xylevel}) it follows that 
all the elements in ${\mathcal B}$ on the same line 
possess the same level.
We let the levels be 
$l^{-1}_t$ (resp. $k^{-1}_s$) on the 
$t$th horizontal (resp. $s$th vertical) line.

The bilinearization of the tropical $R$ attained 
in Section \ref{sec:B-tropR} leads to another formulation 
of the tropical vertex model in terms of tau functions.
In view of Figure \ref{fig:fig2} we
duplicate each line in ${\mathcal L'}$ 
into a pair of parallel lines to form a new lattice ${\mathcal L}$.
Its unit structure looks as Figure \ref{fig:fig2} 
where we now replace $(L,K)$ by $(L(t),K(s))$.
To each face of ${\mathcal L}$ we assign a coordinate 
$(s,t)$ which now takes values in $(\Z/2)^2$ so that 
$\tau^0:(s,t) \in \Z^2$, 
$N:(s,t-\frac{1}{2})$, $W: (s-\frac{1}{2},t)$, 
$S:(s,t+\frac{1}{2})$, $E: (s+\frac{1}{2},t)$, 
$\tau^1:(s+\frac{1}{2},t-\frac{1}{2})$, 
$\tau^2:(s-\frac{1}{2},t-\frac{1}{2})$, 
$\tau^3:(s-\frac{1}{2},t+\frac{1}{2})$, 
$\tau^4:(s+\frac{1}{2},t+\frac{1}{2})$.
There are three types of faces depending on whether 
the coordinate $(s,t)$ belongs to 
$\Z^2$, $(\Z+\frac{1}{2})^2$ or else.
To the face of ${\mathcal L}$ at $(s,t)$,  
we associate a tau function $\tau_i(s,t)$ 
having the components $1 \le i \le n-2$ if $s+t \in \Z+\frac{1}{2}$ and 
$0 \le i \le n$ otherwise.
Then up to a boundary condition,
the tropical vertex model is equivalent to imposing the 
$D^{(1)}_n$ bilinear equations on the tau functions 
around each face at $(s,t) \in \Z^2$.

Let us construct such system of tau functions $\{\tau_i(s,t)\}$ 
by making use of the result in Section \ref{subsec:local},
which may be regarded as a solution of our tropical vertex model.
To each face of ${\mathcal L}$ at $(s,t) \in (\Z/2)^2$ 
we attach the time variables
\begin{align*}
&x_j(s,t) = x(s,t) + z_j,\\
&x(s,t) = \eta +
\sum_{t'\ge t}\veps(L(t')^{-1}) + 
\sum_{t'> t}{\tilde \veps}(L(t')^{-1}) + 
\sum_{s'<s}\veps(K(s')^{-1}) + 
\sum_{s'\le s}{\tilde \veps}(K(s')^{-1}),
\end{align*}
where $z_j$ is defined in (\ref{eq:zdef}) and $\eta = \tilde{\eta}$ 
is an arbitrary odd time.
The sums extend over $s', t' \in \Z$.
Let $g,h \in e^{D_\infty}$, 
$g',h' \in e^{D'_\infty}$ be as in (\ref{eq:mX})--(\ref{eq:mtau0}).
Set
\begin{align*}
\tau_j(s,t)&= F_{1,1}(x_j;g),  \;\;  1 \le j \le n-2, \;\;
s+t \in \Z+1/2,\\
&= \begin{cases}
f_j(x_1;g') & j = 0,1, \\
F_{1,1}(x_j;g) & 2 \le j \le n-2, \\
f_{n-j}(x_1;h') & j = n-1,n,
\end{cases}\quad (s,t) \in (\Z+1/2)^2,\\
&= \begin{cases}
(F_{0,1} + (-1)^jiF_{0,1;1})(x_1;g) & j = 0,1, \\
F_{1,1}(x_{j-1};g) & 2 \le j \le n-2,\\
(F_{0,0} + 
\frac{(-1)^{n-j}ia_{n-1}}{K(s)L(t)}F_{-1,1;1})(x_{n-2};\omega(g))& j = n-1,n,\\
\end{cases}\quad (s,t) \in \Z^2,
\end{align*}
where $x_j = x_j(s,t)$, and the time variable $y$ suppressed here
is taken as explained in the paragraph preceding (\ref{eq:X}).
Since $x(s,t)$ obeys the recursion relation 
corresponding to (\ref{eq:times}), we conclude that 
$\{\tau_j(s,t)\}$ satisfies all 
the $D^{(1)}_n$ bilinear equations if the
parameters $\la, \ka, \alpha, \beta$ at $(s,t) \in \Z^2$ 
are identified with those in Theorem \ref{th:n-finite} 
with $(L,K)$ replaced by $(L(t),K(s))$.
In particular, the levels $l_t^{-1}$ and $k_s^{-1}$ of the edge variables 
are determined by
\begin{equation*}
l_t = L(t)^2\prod_{i=2}^{n-1}(L(t)^2-a^2_i),\qquad 
k_s = K(t)^2\prod_{i=2}^{n-1}(K(s)^2-a^2_i).
\end{equation*}
Note that the formal sums over $s',t'$ in $x(s,t)$ 
make sense under a suitable choice of $\eta$.

\begin{remark}\label{rm:ybe}
Consider the composition of Figure \ref{fig:fig2}
corresponding to $R_1R_2R_1(x,y,z)=(x',y',z')$ and 
$R_2R_1R_2(x,y,z)=(x",y",z")$.

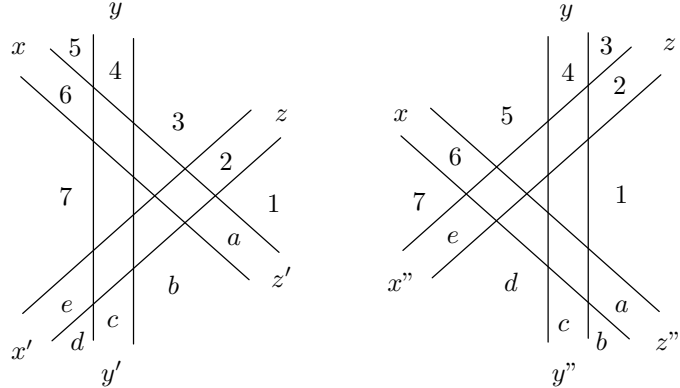
\begin{figure}[h]
\unitlength 0.1in
\begin{picture}( 34.1000, 19.2000)(  9.7000,-28.7000)
%
\special{pn 8}%
\special{pa 1400 1200}%
\special{pa 1400 2800}%
\special{fp}%
\put(23.3000,-25.2000){\makebox(0,0)[lb]{$z'$}}%
\put(14.4000,-30.4000){\makebox(0,0)[lb]{$y'$}}%
\put(9.7000,-29.0000){\makebox(0,0)[lb]{$x'$}}%
\put(12.3000,-26.5000){\makebox(0,0)[lb]{$e$}}%
\put(12.8000,-28.5000){\makebox(0,0)[lb]{$d$}}%
\put(14.7000,-27.3000){\makebox(0,0)[lb]{$c$}}%
\put(17.9000,-25.6000){\makebox(0,0)[lb]{$b$}}%
\put(21.0000,-23.0000){\makebox(0,0)[lb]{$a$}}%
\put(23.5000,-16.5000){\makebox(0,0)[lb]{$z$}}%
\put(14.8000,-11.2000){\makebox(0,0)[lb]{$y$}}%
\put(9.7000,-13.0000){\makebox(0,0)[lb]{$x$}}%
\put(12.2000,-20.9000){\makebox(0,0)[lb]{7}}%
\put(12.2000,-15.6000){\makebox(0,0)[lb]{6}}%
\put(12.7000,-13.1000){\makebox(0,0)[lb]{5}}%
\put(14.8000,-14.3000){\makebox(0,0)[lb]{4}}%
\put(18.1000,-17.0000){\makebox(0,0)[lb]{3}}%
\put(20.6000,-19.1000){\makebox(0,0)[lb]{2}}%
\put(23.1000,-21.3000){\makebox(0,0)[lb]{1}}%
%
\special{pn 8}%
\special{pa 2226 1596}%
\special{pa 1030 2660}%
\special{fp}%
%
\special{pn 8}%
\special{pa 2380 1740}%
\special{pa 1184 2804}%
\special{fp}%
%
\special{pn 8}%
\special{pa 1176 1276}%
\special{pa 2372 2340}%
\special{fp}%
%
\special{pn 8}%
\special{pa 1020 1420}%
\special{pa 2216 2484}%
\special{fp}%
%
\special{pn 8}%
\special{pa 1610 1220}%
\special{pa 1610 2820}%
\special{fp}%
%
\special{pn 8}%
\special{pa 3990 1190}%
\special{pa 3990 2790}%
\special{fp}%
\put(43.4000,-28.5000){\makebox(0,0)[lb]{$z"$}}%
\put(38.0000,-30.4000){\makebox(0,0)[lb]{$y"$}}%
\put(29.4000,-25.4000){\makebox(0,0)[lb]{$x"$}}%
\put(32.5000,-22.9000){\makebox(0,0)[lb]{$e$}}%
\put(35.5000,-25.4000){\makebox(0,0)[lb]{$d$}}%
\put(38.3000,-27.5000){\makebox(0,0)[lb]{$c$}}%
\put(40.3000,-28.5000){\makebox(0,0)[lb]{$b$}}%
\put(41.3000,-26.5000){\makebox(0,0)[lb]{$a$}}%
\put(43.8000,-12.8000){\makebox(0,0)[lb]{$z$}}%
\put(38.4000,-11.3000){\makebox(0,0)[lb]{$y$}}%
\put(29.7000,-16.5000){\makebox(0,0)[lb]{$x$}}%
\put(30.7000,-21.2000){\makebox(0,0)[lb]{7}}%
\put(32.6000,-18.8000){\makebox(0,0)[lb]{6}}%
\put(35.3000,-16.7000){\makebox(0,0)[lb]{5}}%
\put(38.5000,-14.4000){\makebox(0,0)[lb]{4}}%
\put(40.5000,-13.0000){\makebox(0,0)[lb]{3}}%
\put(41.3000,-20.8000){\makebox(0,0)[lb]{1}}%
\put(41.2000,-15.1000){\makebox(0,0)[lb]{2}}%
%
\special{pn 8}%
\special{pa 3166 1586}%
\special{pa 4362 2650}%
\special{fp}%
%
\special{pn 8}%
\special{pa 3010 1730}%
\special{pa 4206 2794}%
\special{fp}%
%
\special{pn 8}%
\special{pa 4216 1266}%
\special{pa 3020 2330}%
\special{fp}%
%
\special{pn 8}%
\special{pa 4370 1410}%
\special{pa 3174 2474}%
\special{fp}%
%
\special{pn 8}%
\special{pa 3780 1210}%
\special{pa 3780 2810}%
\special{fp}%
\end{picture}%
\end{figure}

\vspace{0.2cm}\noindent
According to the bilinearization, assign 
the common time variables and tau functions 
to the faces $1,2,\ldots,7$ of the two diagrams to represent 
the incoming state 
$(x,y,z) \in \B \times \B \times \B$.
Then the time variables on the faces 
$a,b,\cdots,e$ are also the same in the two diagrams.
Therefore one has $(x',y',z') = (x",y",z")$ for the associated 
outgoing state in $\B \times \B \times \B$, 
which implies the Yang-Baxter equation 
$R_1R_2R_1 = R_2R_1R_2$.
\end{remark}

\section{Reduction to $A^{(2)}_{2n-1}$ and $C^{(1)}_n$}\label{sec:red}

Here we introduce the tropical $R$ for $A^{(2)}_{2n-1}$ and $C^{(1)}_n$
based on Remark \ref{rem:RAC} and explain how the
results on the $D^{(1)}_n$ case can be specialized to them.
We note that these tropical $R$ have 
more intrinsic characterization 
in terms of geometric crystals as argued in \cite{KOTY} for $D^{(1)}_n$.

\subsection{\mathversion{bold}
Tropical $R$ for $A^{(2)}_{2n-1}$ and $C^{(1)}_n$}\label{subsec:redR}
In this section we denote the tropical $R$ for $D^{(1)}_n$ by 
$R_n: \B_n \times \B_n \rightarrow \B_n \times \B_n$.
Let $\B' = \{x=(x_1,\ldots, x_n,\xb_n,\ldots, \xb_1) \}$ 
be the set of variables and 
$\rho: \B' \rightarrow \B_{n+1}$ be the embedding 
defined by 
\[
\rho\bigl((x_1,\ldots, x_n,\xb_n,\ldots, \xb_1)\bigr)
=(x_1,\ldots, x_n,1,\xb_n,\ldots, \xb_1).
\]
The tropical $R$ for $A^{(2)}_{2n-1}$ is the birational map 
$R'_n: \B' \times \B' \rightarrow \B' \times \B'$ defined by
\begin{equation*}
R'_n(x,y) = (x',y') \Leftrightarrow 
R_{n+1}(\rho(x),\rho(y)) = (\rho(x'),\rho(y')).
\end{equation*}

Let $\B'' = \{x=(x_0,x_1,\ldots, x_n,\xb_n,\ldots, \xb_1) \}$ 
be the set of variables and 
$\rho': \B'' \rightarrow \B_{n+2}$ be the embedding 
defined by 
\[
\rho'\bigl((x_0,x_1,\ldots, x_n,\xb_n,\ldots, \xb_1)\bigr)
=(x_0,x_1,\ldots, x_n,1,\xb_n,\ldots, \xb_1,x_0).
\]
The tropical $R$ for $C^{(1)}_n$ is the birational map 
$R''_n: \B'' \times \B'' \rightarrow \B'' \times \B''$ defined by
\begin{equation*}
R''_n(x,y) = (x',y') \Leftrightarrow 
R_{n+2}(\rho'(x),\rho'(y)) = (\rho'(x'),\rho'(y')).
\end{equation*}
The maps $R'_n$ and $R''_n$ are well defined due to 
Remark \ref{rem:RAC}.
They satisfy the inversion relation and the Yang-Baxter equation.
We shall write $R'_n$ and $R''_n$ simply as $R'$ and $R''$ 
when the rank needs not be specified.
\subsection{\mathversion{bold} Bilinearization of $R'$ and $R''$}
\label{subsec:RpRpp}

First consider the tropical $R'_{n-1}$ for $A^{(2)}_{2n-3}$.
The bilinear equation for it is given by 
setting $\la_n = \ka_n=1$ and 
specializing the tau functions as 
\begin{equation}\label{eq:taunn-1}
\tau^J_{n-1} = \tau^J_n \qquad 0 \le J \le 4
\end{equation}
in the $D^{(1)}_n$ bilinear equations.
The constraint (\ref{eq:taunn-1}) implies 
$x_n = y_n = x'_n= y'_n = 1$ in 
(\ref{eq:xy-tau})-(\ref{eq:xyprime-tau}), 
which is consistent with Remark \ref{rem:RAC} and 
Proposition \ref{pr:unique-existence}.
It makes the bilinear equations 
$\langle J, n-1 \rangle$ and $\langle J, n \rangle$ 
equivalent.
We call the resulting subsystem 
$\{\langle J,i\rangle \mid 1 \le J \le 4, 0 \le i \le n-1\}$ 
of the $D^{(1)}_n$ case the 
$A^{(2)}_{2n-3}$ bilinear equations.

As for the tropical $R''_{n-2}$ for $C^{(1)}_{n-2}$, 
the bilinear equation is obtained by further setting 
$\la_1=\ol{\la}_1, \ka_1 = \ol{\ka}_1$ and 
imposing the constraint 
\begin{equation}\label{eq:tau01}
\tau^J_{0} = \tau^J_1 \qquad 0 \le J \le 4
\end{equation}
on the $A^{(2)}_{2n-3}$ bilinear equations.
The constraint implies  
$x_1 = \xb_1,\, y_1 = \yb_1$, 
$x'_1 = \xb'_1,\, y'_1 = \yb'_1$ in 
(\ref{eq:xy-tau})-(\ref{eq:xyprime-tau}), which is again
consistent with Remark \ref{rem:RAC} and 
Proposition \ref{pr:unique-existence}.
It makes the equations 
$\langle J, 0 \rangle$ and $\langle J, 1 \rangle$ 
equivalent.
We call the resulting subsystem 
$\{\langle J,i\rangle \mid 1 \le J \le 4, 1 \le i \le n-1\}$ 
of the $A^{(2)}_{2n-3}$ case the 
$C^{(1)}_{n-2}$ bilinear equations.
Under these specializations Theorem \ref{th:bilinear}
is still valid for $R'_{n-1}$ and $R''_{n-2}$.

\subsection{\mathversion{bold}
Solutions of bilinear equations for $R'$ and $R''$}
\label{subsec:solRpRPP}

Results in Section \ref{subsec:local} can be 
specialized to make (\ref{eq:taunn-1}) and (\ref{eq:tau01}) hold.
First we deal with (\ref{eq:taunn-1}) for $R'$.
Let us restrict the elements 
(\ref{eq:gp}) and (\ref{eq:g}) to the case $M=1$ and 
take the limit
\begin{equation}
q'_1 \rightarrow 0, \quad p'_1 \rightarrow a_m
\end{equation}
for some $2 \le m \le n-1$ 
keeping $(1-p'_1/a_m)/{q'_1}^2$ to be an arbitrary 
nonzero constant.
It makes the latter of (\ref{eq:red}) void, leaving 
$a_1,\ldots, a_{n-1}$ as free parameters.
All the tau functions in 
(\ref{eq:mX})-(\ref{eq:mtau0}) are well defined 
in the limit as argued in the end of 
Appendix \ref{subsec:ef}.
Moreover in (\ref{eq:xp}), the coefficients
$p'_1A(p'_1)/q'_1$ and $q'_1/(p'_1A(-p'_1))$ 
in the second sum tend to zero.
Thus the element $h' \in e^{D'_\infty}$ 
in Lemma \ref{lem:red} actually belongs to 
$e^{B'_\infty}$.
Then the property (\ref{eq:f01}) 
ensures that $\tau^J_{n-1} = \tau^J_n$ for $1 \le J \le 4$ in 
(\ref{eq:mtau1234}) in agreement with (\ref{eq:taunn-1}).
The remaining relation $\tau^0_{n-1} = \tau^0_n$ in 
(\ref{eq:mtau0}) follows {}from (\ref{eq:F=0}).

Let us proceed to $R''$, where the further condition 
(\ref{eq:tau01}) should be satisfied.
The reduction {}from $D^{(1)}_n$ case to 
$C^{(1)}_{n-2}$ is done by setting $\forall c_j=0$ 
in (\ref{eq:gp})-(\ref{eq:g}).
Namely we restrict 
$g' \in e^{D'_\infty}$ and $g \in e^{D_\infty}$ 
to $g' \in e^{B'_\infty}$ and $g \in e^{B_\infty}$, 
which reduces the relevant fermion theory 
essentially to the single component one.
The former condition in (\ref{eq:red}) is still present 
although the latter becomes void.
Then obviously $h' \in e^{B'_\infty}$ holds 
in Lemma \ref{lem:red}, therefore 
(\ref{eq:taunn-1}) and (\ref{eq:tau01}) are valid 
due to (\ref{eq:f01}), (\ref{eq:FB}) and (\ref{eq:sym2}).
%

\appendix
\section{$U_i, V_i$ in terms of tau functions}\label{sec:appA}

Let us prove Lemma \ref{lem:UV}.
Our task is to substitute (\ref{eq:xy-tau}) into the functions 
$U_i, V_i$ and simplify the result by means of the bilinear equations 
on the tau functions.
\begin{lemma}\label{lem:easy1}
\begin{equation*}
\begin{split}
&\frac{1}{x_i} + \frac{1}{\yb_i} = 
\frac{\beta\tau^0_1}{\tau^1_0\tau^3_0}\;(i=1),\;\;
\frac{\beta\tau^0_2\tau^2_0\tau^2_1}{N_1W_1}\; (i=2),\;\;
\frac{\beta\tau^0_{i}\tau^2_{i-1}}{N_{i-1}W_{i-1}}
\; (3 \le i \le n-2),\nonumber\\
&\frac{y_n}{x_{n-1}} + \frac{1}{\yb_{n-1}} = 
\frac{\beta\tau^0_{n-1}\tau^1_n\tau^2_{n-2}\tau^2_{n-1}}
{\ka_nN_{n-2}W_{n-2}\tau^1_{n-1}\tau^2_n}.\nonumber
\end{split}
\end{equation*}
\end{lemma}
\begin{proof}
Apply $\langle2,1\rangle,\ldots, \langle2,n-1\rangle$.
\end{proof}
It is straightforward to check 
\begin{lemma}\label{lem:easy2}
For $1 \le i \le n-1$, the ratio $(y_1\cdots y_i)/(x_1\cdots x_{i-1})$ 
is equal to
\begin{align*}
&\frac{\tau^1_1\tau^2_0}{\ka_1N_1}\; \;(i=1),\quad 
\frac{\la_1W_1\tau^2_0\tau^1_2}
{\ka_1\ka_2N_2\tau^1_0\tau^3_0\tau^2_1}\; \;(i=2),  \quad
\frac{\la_1\cdots\la_{i-1} W_{i-1}(\tau^2_0)^2\tau^1_i}
{\ka_1\cdots\ka_iN_i\tau^1_0\tau^3_0\tau^2_{i-1}}\;\; (3 \le i \le n-2), \\
&\frac{\la_1\cdots\la_{n-2} W_{n-2}(\tau^2_0)^2\tau^1_{n-1}}
{\ka_1\cdots\ka_{n-1}\tau^1_0\tau^3_0\tau^2_{n-2}\tau^2_{n-1}}\; \;(i=n-1).
\end{align*}
\end{lemma}
{}From Proposition \ref{pr:involution}, taking 
$\sigma_\ast$ in Lemma \ref{lem:easy2} implies 
\begin{lemma}\label{lem:easy3}
For $1 \le i \le n-2$, the ratio $(\xb_1\cdots\xb_i)/(\yb_1\cdots\yb_{i-1})$ is equal to
\begin{align*}
&\frac{\tau^2_0\tau^3_1}{\lab_1W_1}\;\; (i=1), \quad 
\frac{\kab_1 N_1 \tau^2_0\tau^3_2}{\lab_1\lab_2W_2\tau^1_0\tau^3_0\tau^2_1}
\;\; (i=2), \quad 
\frac{\kab_1\cdots\kab_{i-1} N_{i-1}(\tau^2_0)^2\tau^3_i}
{\lab_1\cdots \lab_i W_i\tau^1_0\tau^3_0\tau^2_{i-1}}\;\; (3 \le i \le n-2).
\end{align*}
In addition one has 
\begin{equation*}
\frac{\xb_1\cdots\xb_{n-1}}{\yb_1\cdots\yb_{n-2}y_n} = 
\frac{\kab_1\cdots\kab_{n-2}\ka_n N_{n-2}(\tau^2_0)^2\tau^1_{n-1}\tau^3_n}
{\lab_1\cdots \lab_{n-1} \tau^1_0\tau^3_0\tau^1_n\tau^2_{n-2}\tau^2_{n-1}}.
\end{equation*}
\end{lemma}

\begin{proof}[Proof of Lemma \ref{lem:UV}] 
First we consider $V_0$ in (\ref{eq:16}).
In view of (\ref{eq:level}) and (\ref{eq:xylevel}), it is expressed as
\begin{equation}\label{eq:V0}
\begin{split}
V_0 &=l^{-1}v_+ +k^{-1}v_-,\\
v_+ &= \sum_{i=1}^{n-2}\frac{y_1\cdots y_i}{x_1\cdots x_{i-1}}
\left(\frac{1}{x_i}+\frac{1}{\yb_i}\right) + 
\frac{y_1\cdots y_{n-1}}{x_1\cdots x_{n-2}}
\left(\frac{y_n}{x_{n-1}} + \frac{1}{\yb_{n-1}}\right),\\
v_- &= \sum_{i=1}^{n-2}\frac{\xb_1\cdots \xb_i}{\yb_1\cdots \yb_{i-1}}
\left(\frac{1}{x_i}+\frac{1}{\yb_i}\right) + 
\frac{\xb_1\cdots\xb_{n-1}}{\yb_1\cdots\yb_{n-2}y_n}
\left(\frac{y_n}{x_{n-1}} + \frac{1}{\yb_{n-1}}\right).
\end{split}
\end{equation}
Using Lemmas \ref{lem:easy1} 
and \ref{lem:easy2} we get
\begin{equation*}
\frac{\tau^1_0\tau^3_0}{\beta\tau^2_0}v_+ = \frac{\tau^0_1\tau^1_1}{\ka_1N_1} + 
\frac{\la_1\tau^2_0\tau^0_2\tau^1_2}
{\ka_1\ka_2N_1N_2} + 
\tau^2_0\sum_{i=3}^{n-2}\frac{\la_1\cdots\la_{i-1}\tau^0_i\tau^1_i}
{\ka_1\cdots\ka_iN_{i-1}N_i}
+\frac{\la_1\cdots\la_{n-2}\tau^2_0\tau^0_{n-1}\tau^1_n}
{\ka_1\cdots\ka_nN_{n-2}\tau^2_n}.
\end{equation*}
Rewrite $\tau^0_i\tau^1_i\,(1 \le i \le n-2)$ and 
$\tau^0_{n-1}\tau^1_n$ appearing here by means of the bilinear equations
$\langle1,1\rangle-\langle1,n-1\rangle$.
Then all but two terms cancel out, leading to
\begin{equation*}
v_+ = \frac{\beta\tau^2_0}{\alpha\tau^1_0\tau^3_0}
\left(\tau^4_0 - \frac{\la_1\cdots\la_n \tau^2_0\tau^4_n}
{\ka_1\cdots\ka_n \tau^2_n}\right).
\end{equation*}
Similarly one can derive 
\begin{equation*}
v_- = \frac{\beta\tau^2_0}{\alpha\tau^1_0\tau^3_0}
\left(-\tau^4_0 + \frac{\kab_1\cdots\kab_{n-1} \tau^2_0\tau^4_n}
{\lab_1\cdots\lab_{n-1} \tau^2_n}\right)
\end{equation*}
{}from Lemmas \ref{lem:easy1}, \ref{lem:easy3} and 
the bilinear equations $\langle3,1\rangle-\langle3,n-1\rangle$.
Substituting these expressions of $v_{\pm}$ into (\ref{eq:V0}), 
we obtain $V_0 = \beta(k-l)\tau^2_0\tau^4_0/(lk\alpha\tau^1_0\tau^3_0)$ 
as desired.
Starting {}from this result on $V_0$, 
one can apply Lemma \ref{lem:easy1} and the bilinear
equations $\langle 3,1\rangle-\langle 3,n-1\rangle$ 
to the recursion (\ref{eq:vind}) successively to 
derive the formulas for $V_1,\ldots, V_{n-1}$.
Then the formulas for 
$V^{\sigma_1}_0$ and $V^{\sigma_\ast}_1,\ldots, 
V^{\sigma_\ast}_{n-1}$  follow 
{}from (\ref{eq:asbefore})-(\ref{eq:*-tau}) and 
Proposition \ref{pr:involution}.
Finally substituting these results 
into (\ref{eq:w}) and using 
$\langle 4,2\rangle-\langle 4,n-2\rangle$, 
one arrives at the desired expression 
for $U_1,\ldots, U_{n-1}$. 
\end{proof}

\section{Tau functions and free fermions}\label{sec:appB}

Here we briefly recall the free fermion approach to 
the theory of tau functions based on \cite{JM}.
Lemma \ref{lem:jm7.6} and Lemma \ref{lem:sapix} play an 
important role in the main text.
\subsection{Single component fermions}\label{subsec:single}
Single component charged fermions 
and the Fock states are given by
\begin{align*}
&[\psi_j,\psi^\ast_k]_+=\delta_{j k},\quad 
[\psi_j,\psi_k]_+=[\psi^\ast_j,\psi^\ast_k]_+=0\quad (j,k \in \Z),\\
&\psi(p) = \sum_j \psi_j p^j,\quad \psi(p)^\ast = \sum_j \psi^\ast_j p^{-j},\\
&\langle l \vert = \langle \hbox{vac} \vert \Psi^\ast_l,\quad 
\vert l \rangle = \Psi_l \vert \hbox{vac} \rangle,\\
&\Psi_l = \begin{cases}
\psi_{l-1}\cdots \psi_1\psi_0, & l >0 \\
1 & l=0\\ \psi^\ast_l \cdots 
\psi^{\ast}_{-2}\psi^{\ast}_{-1} & l < 0
\end{cases},\quad 
\Psi^{\ast}_l = \begin{cases}
\psi^{\ast}_{0}\psi^{\ast}_1\cdots
\psi^{\ast}_{l-1}, & l >0 \\
1 & l=0\\ \psi_{-1} 
\psi_{-2} \cdots \psi_l & l < 0
\end{cases},\\
&\psi_j\vert \hbox{vac} \rangle = 0, \; 
\langle \hbox{vac} \vert \psi_j^\ast = 0\; (j < 0), \quad 
\psi^\ast_j\vert \hbox{vac} \rangle = 0,\; 
\langle \hbox{vac} \vert \psi_j = 0 \; (j\ge 0).
\end{align*}
The neutral fermions are simple combinations thereof:
\begin{align}
&\phi_j = \frac{\psi_j + (-1)^j\psi^\ast_{-j}}{\sqrt{2}},\quad 
{\hat \phi}_j = i\frac{\psi_j - (-1)^j\psi^\ast_{-j}}{\sqrt{2}},
\nonumber\\
&[\phi_j,\phi_k]_+ = [{\hat \phi}_j,{\hat \phi}_k]_+ 
= (-1)^k\delta_{j,-k},\quad 
[\phi_j,{\hat \phi}_k]_+ = 0, 
\nonumber\\
&\phi(p) = \sum_j \phi_j p^j,\quad 
{\hat \phi}(p) = \sum_j {\hat \phi}_j p^j,
\nonumber\\
&\phi_j \vert l \rangle = {\hat \phi}_j \vert l \rangle = 0 
\quad (l=0,1,\; j < 0),
\quad \langle l \vert \phi_j = \langle l \vert {\hat \phi}_j = 0 
\quad (l=0,1,\; j > 0).\label{eq:kesu}
\end{align}
We prepare time variables and their functions:
\begin{align}
&x=(x_1,x_2,x_3,\ldots), \quad {\tilde x} = (x_1,-x_2,x_3,-x_4,\ldots),
\nonumber\\
&\xi(x,p) = \sum_{j\ge 1}x_j p^j,\nonumber\\
&\veps(a) = (a,\frac{a^2}{2},\frac{a^3}{3},\ldots) 
= -{\tilde \veps}(-a). \label{eq:eps}
\end{align}
We say $x$ is {\em odd} if $x = {\tilde x}$.
With the Hamiltonians 
\begin{align*}
&H(x) = \sum_{m \ge 1} \sum_{j \in \Z} x_m \psi_j \psi^\ast_{j+m},\\
&H'(x,y) = \frac{1}{2}\sum_{l=1,3,5,\ldots} \sum_{m \in \Z}
(-1)^{m+1} (x_l\phi_m \phi_{-m-l} + y_l{\hat \phi}_m {\hat \phi}_{-m-l}),
\end{align*}
the fermions exhibit the time evolutions
\begin{align*}
&\hbox{Ad}e^{H(x)} \psi(p) = e^{\xi(x,p)} \psi(p),\; 
\hbox{Ad}e^{H(x)}  \psi^\ast(p) = e^{-\xi(x,p)} \psi^\ast(p),\\
&\hbox{Ad} e^{H'(x,y)} \phi(p) 
= e^{\xi(x,p)} \phi(p),\;
\hbox{Ad}e^{H(x,y)} {\hat \phi}(p) 
= e^{\xi(y,p)} {\hat \phi}(p)\;   \hbox{ for } 
(x,y) = ({\tilde x}, {\tilde y}),
\end{align*}
where $\hbox{Ad}g(\cdot) = g(\cdot) g^{-1}$.

\subsection{2 component fermions}
\begin{align*}
&[\psi^{(\alpha)}_j,\psi^{(\beta)\ast}_k]_+ = 
\delta_{\alpha \beta}\delta_{jk}, 
[\psi^{(\alpha)}_j,\psi^{(\beta)}_k]_+ = 
[\psi^{(\alpha)\ast}_j,\psi^{(\beta)\ast}_k]_+ = 0\; 
(\alpha,\beta=1,2,\; j,k \in \Z),\\
&\psi^{(\alpha)}(p) = \sum_j \psi^{(\alpha)}_jp^j,\;  
\psi^{(\alpha)\ast}(p) = \sum_j \psi^{(\alpha)\ast}_jp^{-j}, \\
&\phi^{(\alpha)}_j = 
\frac{\psi^{(\alpha)}_j + (-1)^j\psi^{(\alpha)\ast}_{-j}}{\sqrt{2}},\quad 
{\hat \phi}^{(\alpha)}_j = i
\frac{\psi^{(\alpha)}_j - (-1)^j\psi^{(\alpha)\ast}_{-j}}{\sqrt{2}},\\
&\phi^{(\alpha)}(p) = \sum_j \phi^{(\alpha)}_j p^j,\;
{\hat \phi}^{(\alpha)}(p) = \sum_j {\hat \phi}^{(\alpha)}_j p^j,\\
&H(x,y) = \sum_{l\ge 1} \sum_j 
(x_l \psi^{(1)}_j\psi^{(1)\ast}_{j+l} + 
y_l \psi^{(2)}_j\psi^{(2)\ast}_{j+l}),\\
&\psi^{(\alpha)}_j\vert \hbox{vac} \rangle = 0, \; 
\langle \hbox{vac} \vert \psi_j^{(\alpha)\ast} = 0\; (j < 0), \quad 
\psi^{(\alpha)\ast}_j\vert \hbox{vac} \rangle = 0,\; 
\langle \hbox{vac} \vert \psi^{(\alpha)}_j = 0 \; (j\ge 0),\\
&\langle l_1, l_2 \vert = \langle \hbox{vac} \vert 
\Psi^{(1)\ast}_{l_1}\Psi^{(2)\ast}_{l_2}
\quad 
\vert l_2,l_1 \rangle = \Psi^{(2)}_{l_2}\Psi^{(1)}_{l_1}
\vert \hbox{vac} \rangle,
\end{align*}
where $\Psi^{(\alpha)}_l$ and $\Psi^{(\alpha)\ast}_l$ stand for 
$\Psi_l$ and $\Psi^\ast_l$ with $\psi_j, \psi^\ast_j$ 
replaced by $\psi^{(\alpha)}_j, \psi^{(\alpha)\ast}_j$.

\subsection{\mathversion{bold}Algebras $A_\infty, 
B'_\infty, B_\infty, D'_\infty$ and $D_\infty$}
The algebras $B'_\infty \subset D'_\infty$ are defined within the 
single component theory as
$B'_\infty = \{ \sum a_{jk}:\phi_j\phi_k: + d  \mid 
{}^{\exists}N, a_{jk} = 0 
\hbox{ if } \vert j+k \vert > N \}$
and 
$D'_\infty = \{ \sum a_{jk}:\phi_j\phi_k: 
+ \sum b_{jk}:{\hat \phi}_j{\hat \phi}_k:  
+ \sum c_{jk}:{\phi}_j{\hat \phi}_k:  + d  \mid 
{}^{\exists}N, a_{jk} = b_{jk} = c_{jk} = 0 
\hbox{ if } \vert j+k \vert > N \}$.
The algebra $A_\infty$ in the 2 component theory is defined as
$A_\infty = 
\{\sum a_{\alpha i, \beta j} :\psi^{(\alpha)}_i\psi^{(\beta)\ast}_j: + d 
\mid {}^\exists N, a_{\alpha i, \beta j} = 0 
\hbox{ if } \vert i-j \vert > N\}$.
The automorphism $\sigma, \kappa, \omega$ of $A_\infty$
and the homomorphism $\iota: D'_\infty \rightarrow A_\infty$ are given by
\begin{align}
&\sigma: \psi^{(\alpha)}_j \mapsto (-1)^j\psi^{(\alpha)\ast}_{-j},\;
\psi^{(\alpha)\ast}_j \mapsto (-1)^j\psi^{(\alpha)}_{-j},\nonumber\\
&\kappa: \psi^{(\alpha)}_j \mapsto i\psi^{(\alpha)}_j,\;
\psi^{(\alpha)\ast}_j \mapsto -i\psi^{(\alpha)\ast}_j,\nonumber\\
&\omega: \psi^{(\alpha)}_j \mapsto \psi^{(\alpha)}_{j-1},\;
\psi^{(\alpha)\ast}_j \mapsto \psi^{(\alpha)\ast}_{j-1},
\label{eq:omega}\\
&\iota: \phi_j \mapsto \phi^{(1)}_j,\; 
{\hat \phi}_j \mapsto \phi^{(2)}_j.\nonumber
\end{align}
The algebras $B_\infty \subset D_\infty$ are defined as
$D_\infty = \{ X \in A_\infty \mid \sigma(X) = X \}$ and 
$B_\infty = \{ X \in D_\infty \mid \pi(X) = X \}$, 
where $\pi$ is the projection to the first component, i.e.,  
$\pi(\psi^{(\alpha)}_j) = \delta_{\alpha 1}\psi^{(\alpha)}_j$, 
$\pi(\psi^{(\alpha)\ast}_j) = \delta_{\alpha 1}\psi^{(\alpha)\ast}_j$.
Thus $B_\infty$ actually stays within the single component theory.
The map 
\begin{equation*}
\iota + \kappa \iota: D'_\infty \longrightarrow D_\infty,
\quad B'_\infty \longrightarrow B_\infty
\end{equation*}
is an isomorphism and one has 
\begin{equation}\label{eq:HHH}
\iota(H'(x,y)) + \kappa\iota(H'(x,y))
= H(x,y)\quad 
\hbox{ if } (x,y) = ({\tilde x}, {\tilde y}).
\end{equation}
The set of well-defined group elements corresponding to these 
algebras will be denoted, by abuse of notation, by
\begin{equation*}
e^{D_\infty} := \{ e^{X_1}\cdots e^{X_k} \mid k \ge 0, \,
X_1,\ldots, X_k \in D_\infty \hbox{ are locally nilpotent} \},
\end{equation*}
and similarly for $A_\infty, B_\infty, B'_\infty$ and $D'_\infty$.
For $g' \in e^{D'_\infty}$, 
we write 
\begin{equation}\label{eq:ggp}
g = e^{\iota(X') + \kappa\iota(X')} 
= \iota(g')\kappa\iota(g') \in e^{D_\infty}.
\end{equation}

\subsection{Tau functions}
We concern two kinds of tau functions:
\begin{equation}\label{eq:fF}
\begin{split}
&f_l(x,y;h) = \langle l \vert 
e^{H'(x,y)}h\vert l \rangle \quad 
\hbox{ for } (x,y) = ({\tilde x}, {\tilde y}),\; 
h \in e^{D'_\infty},\\
&F_{l_1,l_2;l}(x,y;g) = 
\langle l_1,l_2\vert e^{H(x,y)}g\vert l_2-l,l_1+l \rangle\quad 
g \in e^{A_\infty}.
\end{split}
\end{equation}
$F_{l_1,l_2:0}(x,y;g)$ will simply be denoted by 
$F_{l_1,l_2}(x,y;g)$.
They enjoy the symmetry
\begin{align}
&F_{l_1,l_2;l}(x,y;g) = (-1)^{(l_1+l_2)l}
F_{1-l_1,1-l_2;-l}({\tilde x},{\tilde y};\sigma(g)),\label{eq:sym1}\\
&F_{l_1,l_2;l}(x,y_;g) = F_{l_1-1,l_2-1;l}(x,y;\omega(g)).\label{eq:sym2}
\end{align}
Clearly we have 
\begin{equation}\label{eq:FB}
F_{l_1,l_2;l}(x,y;g) = 0 \qquad \hbox{ if } \, l \neq 0 \, 
\hbox{ and } \, g \in B_\infty.
\end{equation}
$F_{l_1,l_2;l}(x,y_;g)$ will be denoted by 
$F_{l_1,l_2;l}(x,y)$ or $F_{l_1,l_2;l}(x)$ for simplicity.

\subsection{Bilinear identities}
For any $g \in e^{A_\infty}$ one has the hierarchy of equations
\begin{small}
\begin{equation*}
\begin{split}
&\oint\frac{dk}{2\pi i k}(-)^{l_2+l_2'}k^{l_1-l_1'-1}
e^{\xi(x-x',k)}F_{l_1-1,l_2; l+1}(x-\veps(k^{-1}),y)
F_{l_1'+1,l_2'; l'-1}(x'+\veps(k^{-1}),y')\\
+&\oint\frac{dk}{2\pi i k}k^{l_2-l_2'-1}
e^{\xi(y-y',k)}F_{l_1,l_2-1; l}(x,y-\veps(k^{-1}))
F_{l_1',l_2'+1; l'}(x',y'+\veps(k^{-1}))=0,
\end{split}
\end{equation*}
\end{small}
\hspace{-0.2cm}
where $l_1-l_1'\ge l'-l \ge l_2'-l_2+2$ and 
the integration is taken along a small contour around $k=\infty$. 
This is eq.(4.4) in \cite{JM}.
Let $b_1, b_2, b_3$ be nonzero constants.
We take $y=y'$ and suppress the dependence on it.
The above equations contain the following:
\begin{small}
\begin{align}
&(b_2^{-1}-b_3^{-1})F_{l_1,l_2;l}(x+\veps(b_2^{-1})+\veps(b_3^{-1}))
F_{l_1-1,l_2+1;l+1}(x+\veps(b_1^{-1})) + 
\hbox{cyc} = 0,\label{eq:bl1}\\
&(b_2^{-1}-b_3^{-1})F_{l_1,l_2;l}(x+\veps(b_2^{-1})+\veps(b_3^{-1}))
F_{l_1-1,l_2;l+1}(x+\veps(b_1^{-1})) + 
\hbox{cyc} = 0,\label{eq:bl2}\\
&(b_2^{-1}-b_3^{-1})F_{l_1,l_2;l}(x+\veps(b_2^{-1})+\veps(b_3^{-1}))
F_{l_1-1,l_2;l}(x+\veps(b_1^{-1})) + 
\hbox{cyc} = 0,\label{eq:bl3}\\
&(b_2-b_3)F_{l_1,l_2;l}(x+\veps(b_2^{-1})+\veps(b_3^{-1}))
F_{l_1,l_2;l}(x+\veps(b_1^{-1})) + 
\hbox{cyc} = 0,\label{eq:bl4}
\end{align}
\end{small}
\hspace{-0.25cm}
where $+ \hbox{cyc}$ means the additional two terms obtained 
by the cyclic permutations of $b_1, b_2, b_3$.
They also contain
\begin{equation}\label{eq:blaux}
\begin{split}
&F_{0,1}(x+\veps(b_2^{-1}))F_{0,0}(x+\veps(b_1^{-1}))-
F_{0,1}(x+\veps(b_1^{-1}))F_{0,0}(x+\veps(b_2^{-1}))\\
&=(b_1^{-1}-b_2^{-1})
F_{1,0;-1}(x+\veps(b_1^{-1})+\veps(b_2^{-1}))
F_{-1,1;1}(x).
\end{split}
\end{equation}
When both $x$ and $y$ are odd, 
$F_{l_1,l_2;l}(x,y;g)$ with $g \in e^{D_\infty}$ 
can be expressed in terms of 
$f_l(x,y;g')$ as noted in eq.(7.6) in \cite{JM}.
We recall this fact in 

\begin{lemma}\label{lem:jm7.6}
Let $g \in e^{D_\infty}$ and $g' \in e^{D'_\infty}$ be related 
as in (\ref{eq:ggp}).
Suppose that $x$ and $y$ are both odd.
Denoting $F_{l_1,l_2;l}(x,y;g)$ and $f_l(x,y;g')$ 
by $F_{l_1,l_2;l}(x)$ and $f_l(x)$, one has 
\begin{align*}
&F_{1,1}(x) = F_{0,0}(x) = f_0(x)f_1(x), \\
&F_{0,1}(x) = F_{1,0}(x) = \frac{1}{2}(f_0(x)^2+f_1(x)^2),\\
&F_{0,1;1}(x) = -F_{1,0;-1}(x) = \frac{i}{2}(f_0(x)^2-f_1(x)^2).
\end{align*}
\end{lemma}
\begin{proof}
The first equality in each line is due to (\ref{eq:sym1}).
Since $\phi_0^2 = {\hat \phi}^2_0 = 1/2$, we may write 
$\hbox{Ad}e^{H'(x,y)}g' =$ $X_0 + X_1\phi_0 + X_2{\hat \phi}_0$  
$+ X_3\phi_0{\hat \phi}_0$,
where $X_j = X_j(\phi,{\hat \phi})$ contains 
$\phi_k, {\hat \phi}_k$ only for $k\neq 0$ 
and its order is even (resp. odd)
for $j=0,3$ (resp. $j=1,2$).
Due to (\ref{eq:kesu}) we find 
$f_l(x;g') = \langle X_0 \rangle + (-1)^li\langle X_3 \rangle/2$
for $l=0,1$, where 
$\langle X_j \rangle  := 
\langle 0 \vert X_j \vert 0 \rangle = 
\langle 1 \vert X_j \vert 1 \rangle$.
On the other hand {}from (\ref{eq:HHH}) and (\ref{eq:ggp}) 
it follows that 
$\hbox{Ad}e^{H(x,y)}g = 
(X'_0 + X'_1\phi^{(1)}_0 + X'_2\phi^{(2)}_0 
+ X'_3\phi^{(1)}_0\phi^{(2)}_0)
({\hat X}'_0 + {\hat X}'_1{\hat \phi}^{(1)}_0 
+ {\hat X}'_2{\hat \phi}^{(2)}_0 
+ {\hat X}'_3{\hat \phi}^{(1)}_0{\hat \phi}^{(2)}_0)$,
where $X'_j$ and ${\hat X}'_j$ are obtained {}from 
$X_j$ by substituting the 2 component fermions 
into $X_j$ as $X'_j = X_j(\phi^{(1)},\phi^{(2)})$, 
${\hat X}'_j = X_j({\hat \phi}^{(1)},{\hat \phi}^{(2)})$.
Due to the property (\ref{eq:kesu}) for each component, 
we get 
$F_{0,0}(x;g) = \langle 0,0 \vert (X'_0 {\hat X}'_0 + 
X'_3 {\hat X}'_3\phi^{(1)}_0\phi^{(2)}_0 
{\hat \phi}^{(1)}_0{\hat \phi}^{(2)}_0) \vert 0,0 \rangle = 
\langle X_0 \rangle^2 + \langle X_3 \rangle^2/4 = 
(\langle X_0 \rangle + i\langle X_3 \rangle/2)
(\langle X_0 \rangle - i\langle X_3 \rangle/2)$, proving the 
first relation.
The other relations can be shown similarly.
\end{proof}

For odd $x$ and $y$, it is clear {}from the above proof that 
\begin{equation}\label{eq:f01}
f_0(x,y;g') = f_1(x,y;g') \qquad \hbox{ if } \, g' \in e^{B'_\infty}.
\end{equation}
Lemma \ref{lem:jm7.6} admits a generalization to 
\begin{lemma}\label{lem:sapix}
Let $g \in e^{D_\infty}$ and $g' \in e^{D'_\infty}$ be 
as in (\ref{eq:ggp}).
Suppose $y$ is odd and 
${\tilde x} = x + \veps(c^{-1}) - {\tilde \veps}(c^{-1})$.
Denoting $F_{l_1,l_2;l}(x,y;g)$ and $f_l(x,y;g')$ 
by $F_{l_1,l_2;l}(x)$ and $f_l(x)$, one has 
\begin{align*}
&\left.\begin{array}{r}
2F_{0,1}(x) \\ 
-2iF_{0,1;1}(x)
\end{array}\right\} = f_0(x+\veps(c^{-1}))f_0(x-{\tilde \veps(c^{-1})}) \pm
f_1(x+\veps(c^{-1}))f_1(x-{\tilde \veps(c^{-1})}),\\
&\left.\begin{array}{r}
2F_{0,0}(x) \\ 
2ic^{-1}F_{-1,1;1}(x)
\end{array}\right\} = f_0(x+\veps(c^{-1}))f_1(x-{\tilde \veps(c^{-1})}) \pm
f_1(x+\veps(c^{-1}))f_0(x-{\tilde \veps(c^{-1})}).
\end{align*}
\end{lemma}
\begin{proof}
In (\ref{eq:bl1})--(\ref{eq:bl3}), set $b_1=-b_2=c$ and $b_3^{-1}=0$.
Writing $x_+=x+\veps(c^{-1}), x_-=x-{\tilde \veps}(c^{-1})$, one has
\begin{align*}
&F_{1,1}(x_-)F_{0,1}(x_+) + F_{1,1}(x_+)F_{0,1}(x_-) 
= 2F_{1,1}({\tilde x})F_{0,1}(x),\\
&F_{1,1}(x_-)F_{0,1;1}(x_+) + F_{1,1}(x_+)F_{0,1;1}(x_-) 
= 2F_{1,1}({\tilde x})F_{0,1;1}(x),\\
&F_{1,0;-1}(x_-)F_{0,1}(x_+) + F_{1,0;-1}(x_+)F_{0,1}(x_-) 
= 2F_{1,0;-1}({\tilde x})F_{0,1}(x).
\end{align*}
By applying (\ref{eq:sym1}) together with 
$\sigma(g) = g$, the right hand sides here
are equal to $2F_{0,0}(x)F_{0,1}(x), 
2F_{0,0}(x)F_{0,1;1}(x)$ and 
$-2F_{0,1;1}(x)F_{0,1}(x)$, respectively.
Since both $x_+$ and $x_-$ are odd, the left hand sides
are expressed in terms of $f_0(x_{\pm}), f_1(x_\pm)$ by
Lemma \ref{lem:jm7.6}.
Solving the resulting three equations with respect to 
$F_{0,0}(x), F_{0,1}(x)$ and $F_{0,1;1}(x)$, we 
obtain the formulas in question up to an overall sign, which can be fixed 
by considering the smooth limit $g \rightarrow 1$.
The formula for $F_{-1,1;1}(x)$ can be derived similarly 
by using (\ref{eq:blaux}).
\end{proof}

Note that in the limit $c \rightarrow \infty$, 
the first three relations in Lemma \ref{lem:sapix} 
reduce to Lemma \ref{lem:jm7.6}.

\subsection{Explicit formula}\label{subsec:ef}
Let us write down the explicit formula for the 
vacuum expectation values (\ref{eq:fF}).
To simplify the result we introduce the parameters 
$p_1,\ldots, p_{2N+M}$ and $q_1,\ldots, q_M$ 
to modify the parameterization of 
$g' \in e^{D'_\infty}$ and $g \in e^{D_\infty}$ in (\ref{eq:gp}) as 
\begin{align*}
&g' = \exp\Bigl(\sum_{i=1}^N b_i\phi(p_i)\phi(p_{\bar i}) 
+ \sum_{j=1}^M c_j\phi(p_{\tilde j}){\hat \phi}(q_j)\Bigr),\\
&g = \exp\Bigl(\sum_{i=1}^Nb_i\bigl(
\psi^{(1)}(p_i)\psi^{(1)\ast}(-p_{\bar i}) -
\psi^{(1)}(p_{\bar i})\psi^{(1)\ast}(-p_i)\bigr)\\
&\qquad + \sum_{j=1}^M c_j\bigl(
\psi^{(1)}(p_{\tilde j})\psi^{(2)\ast}(-q_j) - 
\psi^{(2)}(q_j)\psi^{(1)\ast}(-p_{\tilde j})\bigr)\Bigr),
\end{align*}
where ${\bar i} = 2N+1-i$ and ${\tilde j} = 2N+M+1-j$.
Given $I \subseteq \{1,\ldots, N\}$ and $J \subseteq \{1,\ldots, M\}$,
we shall write ${\bar I} = \{ {\bar i} \mid i \in I\}$ and 
${\tilde J} = \{ {\tilde j} \mid j \in J\}$.
The cardinality of a set $I$ is denoted by $\vert I \vert$.
For any sets $K$ and $K'$ we set 
\begin{align*}
&\Delta^{\pm}_K(p) = 
\prod_{\mu, \nu \in K,\, \mu < \nu}(p_\mu \pm p_\nu),\\
&\Delta_{K,K'}^{m,l}(p) = \frac{\Delta^-_K(p)\Delta^-_{K'}(-p)}
{\prod_{\mu \in K, \nu \in K'}(p_\mu + p_\nu)}
\prod_{\mu\in K}p_\mu^{m+l}
\prod_{\mu\in K'}(-p_\mu)^{-m-l+1},
\end{align*}
and similarly for $q$. Let $x$ and $y$ be odd time variables.
Then we have 
\begin{align*}
&f_l(x,y;g') = \sum_{I,J} \veps_{l,J}2^{-\vert I \vert - \vert J \vert}
b_I(x)c_J(x,y)\frac{\Delta^-_K(p) \Delta^-_J(q)}
{\Delta^+_K(p) \Delta^+_J(q)} \qquad l = 0,1,\\
&K = I \sqcup {\bar I} \sqcup {\tilde J},\qquad 
\veps_{l,J} = \begin{cases}
1 & \vert J \vert \hbox{ even}\\
i(-1)^l & \vert J \vert \hbox{ odd},
\end{cases}\\
&b_I(x) = \prod_{i \in I} b_i\exp(\xi(x,p_i) - \xi(x,-p_{\bar i})),\\
&c_J(x,y) = \prod_{j \in J}c_j
\exp(\xi(x,p_{\tilde j}) + \xi(y,q_j)),
\end{align*}
where the sum $\sum_{I,J}$ extends over all the subsets 
$I \subseteq \{1,\ldots, N\}$ and 
$J \subseteq \{1,\ldots, M\}$.
For any integers $l_1,l_2$ and $l$, 
$F_{l_1,l_2;l}(x,y;g)$ with odd $y$ 
($x$ is not restricted to be odd) reads
\begin{align*}
&F_{l_1,l_2;l}(x,y;g)\\
&= \sum_{I,I',J,J'}(-1)^{\vert I' \vert 
+ (\vert I \vert + \vert I' \vert)\vert J \vert}
b_I(x)b'_{I'}(x)c_J(x,y)c'_{J'}(x,y)
\Delta_{K,K'}^{l_1,l}(p)\Delta_{J',J}^{l_2,-l}(q),\\
&K = I \sqcup {\bar I}' \sqcup {\tilde  J},\qquad 
K' = I' \sqcup {\bar I} \sqcup {\tilde  J}',\\
&b'_{I'}(x) = \prod_{i \in I'} b_i\exp(\xi(x,p_{\bar i}) - \xi(x,-p_i)),\\
&c'_{J'}(x,y) = \prod_{j \in J'}c_j
\exp(-\xi(x,-p_{\tilde j}) + \xi(y,q_j)),
\end{align*}
where the sum $\sum_{I,I',J,J'}$ extends over 
all the subsets 
$I, I' \subseteq \{1,\ldots, N\}$ and 
$J, J' \subseteq \{1,\ldots, M\}$ 
such that $\vert J'\vert - \vert J \vert = l$.

Consider the dependence of $f_l(x,y;g')\, (l=0,1)$ 
and $F_{l_1,l_2;l}(x,y;g)$ 
on $q_1$ when $M=1$.
Apart {}from the factors $c_J(x,y)$ and $c'_{J'}(x,y)$, 
$f_l$ is independent of $q_1$ and 
$F_{l_1,l_2;l}$ is 
proportional to the single power
$q_1^{l_2-l}$.
Therefore all the functions $f_l$ and $F_{l_1,l_2;l}$ 
appearing in (\ref{eq:mX})-(\ref{eq:mtau0}) are 
well defined in the limit $q_1 \rightarrow 0$ with $M=1$.
In particular {}from (\ref{eq:sym2}) one has 
\begin{equation}\label{eq:F=0}
\lim_{q_1 \rightarrow 0}F_{-1,1;1}(x,y;\omega(g)) = 
\lim_{q_1 \rightarrow 0}F_{0,2;1}(x,y;g) = 0.
\end{equation}
These facts are used in Section \ref{subsec:solRpRPP}.

\section*{Acknowledgements}
The authors thank R. Inoue for a careful reading of the manuscript.
A.K. thanks M. Jimbo for explaining eq.(7.6) in \cite{JM}. 
M.O. and Y.Y. were partially supported by 
Grant-in-Aid for Scientific Research JSPS No.14540026 and 
No.11440047, respectively.

\end{document}